\documentclass[11pt]{amsart}

\usepackage{amsmath,amssymb,amsthm}
\usepackage{amsxtra,amsfonts,amscd,euscript}
\usepackage{mathrsfs}

\usepackage{graphicx,epstopdf}
\usepackage{url}
\usepackage{hyperref}
\newtheorem{theorem}{Theorem}
\newtheorem{lemma}{Lemma}
\numberwithin{lemma}{section}
\newtheorem{proposition}[lemma]{Proposition}

\theoremstyle{definition}

\theoremstyle{remark}
\newtheorem{remark}{Remark}

\numberwithin{equation}{section}
\providecommand{\abs}[1]{\left\lvert#1\right\rvert}
\providecommand{\tabs}[1]{\lvert#1\rvert}
\providecommand{\norm}[1]{\left\lVert#1\right\rVert}
\DeclareMathOperator{\re}{\mathrm{Re}}
\DeclareMathOperator{\im}{\mathrm{Im}}

\newcommand{\Z}{\ensuremath{\mathbb{Z}}}

\newcommand{\R}{\ensuremath{\mathbb{R}}}
\newcommand{\C}{\ensuremath{\mathbb{C}}}

\renewcommand{\(}{\left(}
\renewcommand{\)}{\right)}
\renewcommand{\{}{\left\lbrace}
\renewcommand{\}}{\right\rbrace}
\newcommand{\vep}{\varepsilon}
\newcommand{\ud}{\mathrm{d}}
\newcommand{\vp}{\varphi}
\newcommand{\bk}{\backslash}

\DeclareMathOperator{\littleo}{o}
\DeclareMathOperator{\bigO}{O}

\DeclareMathOperator{\schwarz}{\mathcal{S}}

\newcommand{\U}{\mathbb{H}}

\newcommand{\moduli}{\mathfrak{M}}
\newcommand{\schottky}{\mathfrak{S}}
\newcommand{\teich}{\mathfrak{T}}
\newcommand{\univcurve}{\mathscr{M}}
\newcommand{\univschottky}{\mathscr{S}}

\newcommand{\diffform}{\mathcal{E}}
\newcommand{\autoform}{\mathcal{A}}
\newcommand{\holoautoform}{\mathcal{H}}
\newcommand{\Ltwoautoform}{\mathfrak{H}}

\newcommand{\canclass}{\omega}

\newcommand{\zbar}{\bar{z}}

\newcommand{\wbar}{\bar{w}}

\newcommand{\del}{\partial}
\newcommand{\delbar}{\bar{\partial}}

\newcommand{\delbarprime}{\bar{\partial} '}
\newcommand{\var}{\delta}
\newcommand{\varbar}{\bar{\delta}}

\providecommand{\define}[1]{\emph{#1}}

\begin{document}

\title[Holomorphic factorization of determinants]
{Holomorphic factorization of determinants of
Laplacians on Riemann surfaces and a higher genus
generalization of Kronecker's first limit formula}
%Short title in []

\author{Andrew McIntyre}
\address{Centre de Recherches Math\'ematiques, Universit\'e de Montr\'eal,
 C.~P. 6128 Centre-ville station, Montr\'eal QC, H2X 2P1 Canada }
\email{mcintyre@crm.umontreal.ca}
\urladdr{http://www.crm.umontreal.ca/\textasciitilde mcintyre/}
\thanks{}

\author{Leon A. Takhtajan}
\address{Department of Mathematics, Stony Brook University, Stony Brook NY, 11794-3651, USA}
\email{leontak@math.sunysb.edu}
\urladdr{http://www.math.sunysb.edu/\textasciitilde leontak/}
\thanks{}

\date{October 12, 2004}
%Replace with explicit date if necessary
\dedicatory{}

\begin{abstract} 
For a family of compact Riemann surfaces
$X_{t}$ of genus $g>1$, parameterized by
the Schottky space $\mathfrak{S}_{g}$, we
define a natural basis of
$H^{0}(X_{t},\canclass^{n}_{X_{t}})$ which varies
holomorphically with $t$ and generalizes the
basis of normalized abelian differentials of 
the first kind for $n=1$. We introduce a
holomorphic function $F(n)$ on $\schottky_{g}$
which generalizes the classical product
$\prod_{m=1}^{\infty}(1-q^{m})^{2}$ for $n=1$
and $g=1$. We prove the holomorphic
factorization formula
\begin{equation*}
\frac{{\det}'\Delta_{n}}{\det N_{n}}=
 c_{g,n}
 \exp\{-\frac{6n^{2}-6n+1}{12\pi}S\}
  \abs{F(n)}^{2},
\end{equation*}
where 
${\det}'\Delta_{n}$ is the zeta-function
regularized determinant of the Laplace 
operator $\Delta_{n}$ in the hyperbolic metric 
acting on $n$-differentials, $N_{n}$ is the Gram 
matrix of the natural basis with respect to the inner 
product given by the hyperbolic metric, $S$ is 
the classical Liouville action --- a K\"{a}hler 
potential of the Weil-Petersson metric on 
$\schottky_{g}$ --- and  $c_{g,n}$ is a constant 
depending only on $g$ and $n$. The factorization 
formula reduces to Kronecker's first limit formula 
when $n=1$ and $g=1$, and to Zograf's 
factorization formula for $n=1$ and $g>1$.   
\end{abstract}

\maketitle

%----------------------------------------
%   Body
%----------------------------------------

%----------------------------------------
\section{Introduction}
%----------------------------------------

Let $s$ and $\tau$ be complex numbers with $\re s>1$, 
$\im\tau>0$, and define
\[
E(\tau,s)=\sum_{\substack{(m,n)\in\Z^{2}\\(m,n)\neq (0,0)}}
\frac{(\im\tau)^{s}}{\abs{m+n\tau}^{2s}}.
\]
This series was introduced by Kronecker in 1863; 
see \cite{Weil76:Ellipticfns}. It admits meromorphic 
continuation to the entire $s$-plane with a single 
simple pole at $s=1$, and satisfies the functional 
equation
\begin{equation} \label{functional-equation}
\pi^{-s}\Gamma(s)E(\tau,s)=\pi^{s-1}\Gamma(1-s)E(\tau,1-s),
\end{equation}
where $\Gamma(s)$ is Euler's gamma-function.
Kronecker's first limit formula asserts that
\begin{equation} \label{kronecker-1}
E(\tau,s)=\frac{\pi}{s-1}
 -\pi\log\{\frac{4\im\tau\abs{\eta(\tau)}^{4}}
   {\exp(2\Gamma'(1))}\}  + \bigO(s-1)
\end{equation}
as $s\rightarrow 1$, where $\eta(\tau)$ is the
Dedekind eta-function:
\[
\eta(\tau)=q^{\tfrac{1}{24}}
 \prod_{m=1}^{\infty}(1-q^{m}),\quad q=e^{2\pi i\tau}.
\]
See \cite{Weil76:Ellipticfns} and \cite{Lang87:Ellipticfns}
for the proof, and for applications to number theory.
Equation \eqref{kronecker-1} admits an interpretation in
terms of the spectral geometry of the
elliptic curve $E_{\tau}\simeq L\bk\C$,
$L=\Z+\Z\tau$, which goes back to \cite{RS73:Complextorsion}.
Namely, assign to $E_{\tau}$ the flat metric 
$\displaystyle{\frac{1}{\im\tau}\abs{\ud z}^{2}}$, in 
which the area of $E_{\tau}$ is $1$.  Let 
\[
\Delta_{0}(\tau)=-\im\tau\frac{\del^{2}}{\del z\del\zbar}
\]
be the Laplace operator in this metric on $E_{\tau}$, 
acting on functions.  Its eigenvalues are 
\[
\lambda_{\ell}
=\frac{\pi^{2}\abs{\ell}^{2}}{\im\tau},\quad \ell\in L.
\] 
Its determinant is defined by zeta function regularization: 
the function
$\zeta(\tau,s)=\sum_{\lambda_{\ell}\neq 0}
\lambda_{\ell}^{-s}$, defined initially for 
$\re s>1$, admits meromorphic continuation to 
the entire $s$-plane, and one defines
\[
{\det}^{\prime}\Delta_{0}(\tau)
=\exp\{\left.-\frac{\partial}{\partial s}\right|_{s=0}
 \zeta(\tau,s)\},
\]
where the prime indicates omission of zero eigenvalues.
Since $\zeta(\tau,s)=\pi^{-2s}E(\tau,s)$, 
it follows from \eqref{functional-equation}
and \eqref{kronecker-1} that
\begin{equation}\label{torus-0}
{\det}'\Delta_{0}(\tau)
=4\im\tau\abs{\eta(\tau)}^{4}.
\end{equation}
This formula has been used in string
theory for the one-loop computation in the
perturbative approach of Polyakov (see, e.g., 
\cite{DHoker:Strings} and references therein).
 
We restate \eqref{torus-0} in a form convenient 
for generalization to higher genus. Consider the 
Schottky uniformization of the elliptic
curve: $E_{\tau}\simeq\Gamma\bk\C^{\ast}$, 
where $\Gamma$ is the cyclic group generated 
by the dilation $w\mapsto qw$, with fundamental 
region $D=\{w\in\C^{\ast}:\tabs{q}< \tabs{w}\leq 1\}$.
The push-forward of the Euclidean metric 
$(\im\tau)^{-1}\tabs{\ud z}^{2}$ by the map $w=e^{2\pi iz}$ takes 
the form $\rho(w)\tabs{\ud w}^{2}$, where 
$\rho(w) =(4\pi^{2}\im\tau\tabs{w}^{2})^{-1}$.  Setting
\[
S(\tau)=\frac{i}{2}\iint\limits_{D}
 \abs{\frac{\del\log\rho}{\del w}}^{2}
 \,\ud w\wedge \ud\wbar=4\pi^{2}\im\tau,
\]
we can rewrite \eqref{torus-0} as
\begin{equation} \label{torus-1}
\frac{{\det}'\Delta_{0}(\tau)}{\im\tau}
=4\exp\{-\frac{1}{12\pi}S(\tau)\}\abs{F(q)}^{2},
\end{equation}
where
\begin{equation} \label{torus-F}
F(q)=\prod_{m=1}^{\infty}(1-q^{m})^{2}.
\end{equation}
Note that ${\det}'\Delta_{0}(\tau)$
depends only on the isomorphism class
of $E_{\tau}$,  which in turn depends
only on $q$, and that $\im\tau$ also depends
only on $q$.  Hence 
\eqref{torus-1} is an equality of functions
on $\{q\in\C:0<\tabs{q}< 1\}$.

In this paper we extend
\eqref{torus-1} and \eqref{torus-F} from
elliptic curves to
compact Riemann surfaces of genus
$g>1$, and from functions to
$n$-differentials (sections of the
$n$-th power of the canonical bundle). 
To formulate the main result, which may 
be interpreted as a higher genus generalization
of Kronecker's first limit formula, we first recall
some basic facts about uniformization
of Riemann surfaces and about
Teichm\"{u}ller and Schottky spaces
(see Section \ref{basics} for more detail). 
Each compact Riemann surface $X$
of genus $g>1$ carries a unique hyperbolic
metric (a Hermitian metric of constant
negative curvature $-1$), with respect
to which one can define the Laplace operator
$\Delta_{0}(X)$ acting on functions
on $X$, its zeta function (analogous to 
$\zeta(\tau,s)$ defined above), and its regularized
determinant
${\det}'\Delta_{0}(X)$. The Riemann moduli
space is the set $\moduli_{g}$ of isomorphism
classes of compact Riemann surfaces of genus
$g>1$; it carries a natural structure of a
complex orbifold of dimension $3g-3$.  This
generalizes the space 
$\mathrm{PSL}(2,\Z)\bk\{\tau\in\C:\im\tau>0\}$
of isomorphism classes of elliptic curves.  The
determinant
${\det}'\Delta_{0}$ is a real-analytic function
on $\moduli_{g}$.

Now suppose that the Riemann surface
$X$ is marked, i.e., has a distinguished
canonical system of generators
$\alpha_{1},\dots,\alpha_{g},
\beta_{1},\dots,\beta_{g}$
of the fundamental group
$\pi_{1}(X,x_{0})$, $x_{0}\in X$. With
respect to this marking we may define a
normalized basis
$\varphi_{1},\dotsc,\varphi_{g}$
of the space of holomorphic $1$-forms 
--- abelian differentials of the first kind --- 
by the requirement
$\int_{\alpha_{k}}\varphi_{j}=\delta_{jk}$;
then the period matrix $\boldsymbol\tau$ is defined
by $\boldsymbol\tau_{jk}=\int_{\beta_{k}}\varphi_{j}$.
It satisfies 
$\im\boldsymbol\tau_{jk}
=\langle\varphi_{j},\varphi_{k}\rangle
=\tfrac{i}{2}\int_{X}\varphi_{j}\wedge\overline{\varphi}_{k}$ 
by the Riemann bilinear relations.
The Teichm\"{u}ller space $\teich_{g}$ is
the set of isomorphism classes of marked
Riemann surfaces of genus $g$; it is
the universal cover of $\moduli_{g}$, and it
carries a natural structure of a complex
manifold of dimension $3g-3$ with respect
to which the entries of $\boldsymbol\tau$ are 
holomorphic functions.  For $g>1$, the Teichm\"uller space 
generalizes the upper half-plane 
$\{\tau\in\C:\im\tau>0\}$, and 
$\det\im\boldsymbol\tau$ will play the role of 
the factor $\im\tau$ appearing in 
\eqref{torus-1}.

In fact, $\det\im\boldsymbol\tau$ is a 
well defined function
on the Schottky space $\schottky_{g}$, which is
an intermediate cover of $\moduli_{g}$
($\teich_{g}\to\schottky_{g}\to\moduli_{g}$)
defined as follows.  A marked Schottky group is 
a discrete subgroup $\Gamma$ of
the group of linear fractional transformations
$\mathrm{PSL}(2,\C)$,
with distinguished free generators
$L_{1},\dotsc,L_{g}$ satisfying the following
condition:  there exist $2g$ smooth
Jordan curves $C_{r}$, $r=\pm 1,\dotsc,\pm g$,
which form the oriented boundary of a 
domain $D\subset\widehat\C=\C\cup\{\infty\}$, 
such that
$L_{r}C_{r}=-C_{-r}$, $r=1,\dots,g$.  If 
$\Omega$ is the union of images of 
$D$ under $\Gamma$, then 
$\Gamma\backslash\Omega$ is a compact 
Riemann surface of genus $g$.
According to the classical retrosection 
theorem, every compact Riemann surface
may be realized in this manner; if it is 
marked, the condition $C_{k}$ homotopic
to $\alpha_{k}$ for each $k>0$ 
fixes the marked group up to overall conjugation
in $\mathrm{PSL}(2,\C)$.  The overall conjugation
may be fixed by a normalization condition ---
see section \ref{kleinian-groups}.
The Schottky space $\schottky_{g}$ is 
the space of marked normalized Schottky groups
with $g$ generators.
It is a complex manifold of dimension $3g-3$,
covering $\moduli_{g}$ and with universal
cover $\teich_{g}$, and $\det\im\boldsymbol\tau$ 
is a well defined function on it \cite{Zograf89:FPaper}.  
The Schottky space
$\schottky_{g}$ generalizes
the space $\{q\in\C:0<\tabs{q}<1\}$
discussed above.

Like the Teichm\"{u}ller
space $\teich_{g}$, the Schottky
space $\schottky_{g}$ carries a natural
K\"{a}hler metric, the
Weil-Petersson metric. Its global K\"{a}hler
potential can be explicitly constructed
as follows. Let
$\rho(z)\tabs{\ud z}^{2}$
be the hyperbolic metric on $\Omega$ ---
the pull-back of the hyperbolic metric on 
$X\simeq\Gamma\bk\Omega$. 
Following \cite{TakZog87:Liouvilleaction}, set
\begin{equation}\label{liouville-action}
\begin{split}
S &= \frac{i}{2}
\iint\limits_{D}
  \(\abs{\frac{\del\log\rho}{\del z}}^{2}+\rho\)
  \,\ud z\wedge \ud\zbar 
  \\
&\quad + \frac{i}{2}\sum_{k=2}^{g}
 \oint_{C_{k}}\(\log\rho-\frac{1}{2}\log\abs{L_{k}'}^{2}\)
  \(\frac{L_{k}''}{L_{k}'}\,\ud z -
   \frac{\overline{L_{k}''}}{\overline{L_{k}'}}\,\ud\zbar\)
   \\
&\quad +4\pi\sum_{k=2}^{g}\log\abs{c\(L_{k}\)}^{2},
\end{split}
\end{equation}
where for $\gamma=
\bigl(\begin{smallmatrix}
a&b\\c&d\\
\end{smallmatrix}\bigr)\in\Gamma$,
we denote $c(\gamma)=c$. The function
$S:\schottky_{g}\to\R$ is called the 
\define{classical Liouville action} (see 
\cite{TakZog87:Liouvilleaction}
and \cite{TakTeo03:Liouvilleaction}
for details and motivation).  According
to \cite{TakZog87:Liouvilleaction}, the
function $-S$ is a
K\"{a}hler potential of the Weil-Petersson
metric on $\mathfrak{S}_{g}$, i.e.,
\begin{equation} \label{kahler-potential}
\del\delbar S
=2i\omega_{\scriptscriptstyle{WP}},
\end{equation}
where $\del$ and
$\delbar$ are, respectively, the
$(1,0)$ and $(0,1)$ components of the deRham
differential $\ud$ on $\schottky_{g}$, and
$\omega_{\scriptscriptstyle{WP}}$ is the
symplectic form of the Weil-Petersson metric. 
For $g>1$, the function $S$ on $\schottky_{g}$ will play
the role of the 
function $S(\tau)=-2\pi\log\tabs{q}$ on 
$\{q\in\C:0<\tabs{q}<1\}$ appearing in \eqref{torus-1}.

Now we can formulate the following remarkable 
generalization of \eqref{torus-1} and \eqref{torus-F} 
to higher genus Riemann surfaces.
\begin{theorem}[P. Zograf]\label{zograf-thm}
Let $g>1$, and let
${\det}'\Delta_{0}$,
$\im\boldsymbol\tau$ and $S$ be the functions
on the Schottky space $\schottky_{g}$ 
defined above. Then there exists a 
holomorphic function $F:\schottky_{g}\to\C$
such that
\begin{equation} \label{zograf-formula}
\frac{{\det}'\Delta_{0}}{\det\im\boldsymbol\tau}
=c_{g}\exp\{-\frac{1}{12\pi}S\}\abs{F}^{2},
\end{equation}
where $c_{g}$ is a
constant depending only on $g$. 
For points in
$\mathfrak{S}_{g}$ corresponding to 
Schottky groups $\Gamma$
with exponent of convergence $\delta<1$,
the function $F$ is given by the following absolutely 
convergent product:
\begin{equation} \label{F-0}
F=\prod_{\{\gamma\}}\prod_{m=0}^{\infty}
	\(1-q_{\gamma}^{1+m}\),
\end{equation}
where $q_{\gamma}$ is the multiplier of 
$\gamma\in\Gamma$,
and $\{\gamma\}$ runs over all distinct
primitive conjugacy classes in $\Gamma$,
excluding the identity.
\end{theorem}
See section \ref{kleinian-groups} for the definition
of $\delta$, $q_{\gamma}$, and primitive 
$\gamma$.
The factorization formula \eqref{zograf-formula} 
was proved in \cite{Zograf89:FPaper}, and 
the representation \eqref{F-0} was 
discovered later \cite{Zograf97:FPreprint}.
We will refer to \eqref{zograf-formula} together
with \eqref{F-0} as the
\define{Zograf factorization formula}, or
simply Zograf's formula.
Note that when $g=1$, the theorem still
holds provided that $\Delta_{0}$ and $S$ 
are defined as in the discussion of elliptic curves 
above.
In this case, \eqref{zograf-formula} becomes 
\eqref{torus-1}, and the function $F$ reduces to the 
classical product $\eqref{torus-F}$.

Associated to the Riemann surface
$X$ is the Selberg zeta function
\begin{equation}\label{selberg-zeta}
Z(s)=\prod_{\{\gamma\}}\prod_{m=0}^{\infty}
 \(1-q_{\gamma}^{s+m}\),
\end{equation}
where $\{\gamma\}$ runs over all 
distinct nontrivial primitive
conjugacy classes in a \emph{Fuchsian} group
uniformizing $X$.  Defined initially for $\re s>1$,
the Selberg zeta function admits analytic continuation
to the entire $s$-plane, and, according to 
\cite{DP86:DetLaplacian} and \cite{Sarnak87:DetLaplacian},
\[
{\det}'\Delta_{0}=e^{c_{0}(2g-2)}Z'(1)
\]
for some constant $c_{0}$.
Hence Zograf's formula gives a factorization of $Z'(1)$, 
considered as a function on $\schottky_{g}$. 

To motivate the extension from functions to $n$-differentials
on $X$, we first describe a geometric interpretation of Zograf's formula,
in the context of the Quillen metric and the local index 
theorem for families.  We write $\canclass_{X}$ for the holomorphic
cotangent bundle of $X$, and call a smooth section of 
$\canclass_{X}^{n}$ an $n$-differential.
Let $\univcurve_{g}=\moduli_{g,1}$ 
be the universal curve --- the moduli space of compact Riemann 
surfaces of genus $g>1$ with one marked point --- and let
$p:\univcurve_{g}\to\moduli_{g}$ be the corresponding forgetful 
map.  Denote by $T_{V}\univcurve_{g}$ the vertical holomorphic 
tangent bundle of the fibration $p$, and for each 
positive integer $n$, denote
by $\Lambda_{n}$ the direct image bundle 
$p_{\ast}(T_{V}\univcurve_{g}^{-n})$ over $\moduli_{g}$. 
Then the fibre of
$\Lambda_{n}$ over a point $t\in\moduli_{g}$ is isomorphic to
the vector space $H^{0}(X_{t},\canclass_{X_{t}}^{n})$ of
holomorphic $n$-differentials on the Riemann surface $X_{t}=p^{-1}(t)$.
Let $\lambda_{n}=\det\Lambda_{n}$ be the corresponding determinant 
line bundle over $\moduli_{g}$. The hyperbolic metric on the fibres of 
$p$ defines a natural Hermitian metric on $\Lambda_{n}$ and on hence on
$\lambda_{n}$.  The Quillen metric \cite{Quillen85:CRops} on $\lambda_{n}$
is defined by
\begin{equation*}
\norm{\vp}^{2}_{Q,n}
 =\frac{\norm{\vp}_{n}^{2}}{{\det}'\Delta_{n}}
 =\frac{\det N_{n}}{{\det}'\Delta_{n}},
\end{equation*} 
where $\norm{\,\cdot\,}_{n}$ is the Hermitian metric mentioned above,
$\vp=\vp_{1}\wedge\dots\wedge\vp_{d_{n}}$ is a local holomorphic section
of $\lambda_{n}$ at $t\in\moduli_{g}$, 
$[N_{n}]_{jk}=\langle\vp_{j},\vp_{k}\rangle$ is the Gram matrix of the basis
$\vp_{1},\dots,\vp_{d_{n}}$ of $H^{0}(X_{t},\canclass_{X_{t}}^{n})$, and
$\Delta_{n}$ is the Laplace operator in the hyperbolic metric on $X_{t}$ acting
on $n$-differentials. The Quillen metric has the remarkable property 
that the Chern form of the Hermitian line bundle 
$(\lambda_{n},\norm{\,\cdot\,}_{Q,n})$
over $\moduli_{g}$ is proportional to the Weil-Petersson symplectic form
$\omega_{\scriptscriptstyle{WP}}$:
\begin{equation} \label{local-index}
\delbar\del
\log\frac{\det N_{n}}{{\det}'\Delta_{n}}
  =\frac{6n^{2}-6n+1}{6\pi i}
   \,\omega_{\scriptscriptstyle{WP}}.
\end{equation}
This is the local index theorem for families
(see \cite{BK86:HoloAn,BostJol86:Localindexthm,
TakZog87:Localindexthm}).

Theorem \ref{zograf-thm} together with
\eqref{kahler-potential} constitute a refinement
of \eqref{local-index} in the case $n=1$. Let 
$\vp=\vp_{1}\wedge\dots\wedge\vp_{g}$
be the local holomorphic section
of $\lambda_{1}$ determined by the normalized basis 
$\vp_{1},\dots,\vp_{g}$ 
of abelian differentials of the first kind on $X_{t}$.
Then Theorem \ref{zograf-thm} provides (by means
of the function $F$)
an isometry between the line bundle 
$\lambda_{1}$ with the Quillen metric, and the 
line bundle over $\moduli_{g}$ canonically determined
by carrying the Hermitian metric $\exp\{\tfrac{1}{12\pi}S\}$
(see Section 3 in \cite{Zograf89:FPaper} for details).  
(We have used the fact that ${\det}'\Delta_{n}={\det}'\Delta_{1-n}$,
see e.g.
\cite{TakZog87:Localindexthm}.)  Expressed differently,
Zograf's factorization formula is a ``$\del\delbar$
antiderivative'' of \eqref{local-index}.

Based on \eqref{local-index}, it is natural to expect an analogue 
of Theorem \ref{zograf-thm} to hold for all positive integer $n$. 
However, there are two principal differences between the cases 
$n=1$ and $n>1$.

First, for $n=1$ there is a canonical choice of a lattice of maximal 
rank in $H^{0}(X,\canclass_{X})$ provided by the dual to $H_{1}(X,\mathbb{Z})$,
which gives rise to the classical normalized basis of abelian differentials 
described above.  Topology does not fix such a lattice in 
$H^{0}(X,\canclass_{X}^{n})$ when $n>1$. Nevertheless, using Schottky 
uniformization and corresponding Eichler cohomology groups, we 
construct a natural basis of $H^{0}(X_{t},\canclass^{n}_{X_{t}})$ which 
is canonical up to a choice of basis in a space of polynomials,
varies holomorphically with $t\in\schottky_{g}$, and reduces
to the 
classical normalized basis of abelian differentials of the first kind
when $n=1$.

Second, for $n=1$ the holomorphic quadratic differential on $X=X_{t}$ 
which corresponds to the $(1,0)$ form $\del\log{\det}'\Delta_{0}$ 
at $t\in\schottky_{g}$ is given by a local expression in terms of the 
Green's function of $\delbar_{1}$. However, for $n>1$ the corresponding
local expression is not holomorphic, and a holomorphic projection
must be applied to obtain $\del\log{\det}'\Delta_{n}$,
which makes the entire expression non-local. Still, we prove 
that up to a known ``holomorphic anomaly'', 
(which gives rise to the factor involving the classical 
Liouville action $S$),
$\del\log{\det}'\Delta_{n}$ is given by applying the projection operator 
to
\[
T_{n}(z)=\lim_{z'\to z}\(n\frac{\del}{\del z'} -(1-n)\frac{\del}{\del z}\)
 \(K_{n}(z,z')-\frac{1}{\pi}\frac{1}{z-z'}\),
 \quad z\in\Omega,
\]
where $K_{n}$ is the Green's
function for the $\delbar_{n}$-operator. The advantage of this 
representation is that, although $T_{n}$ fails to be holomorphic,
$\displaystyle{\frac{\del T_{n}}{\del\zbar}}$ can be characterized explicitly, and 
the projection can be avoided by means of a contour integration. 
In this we make rigorous the heuristic outline given in 
\cite{Mar87:CFTonRS}, where
$T_{n}$ arises as the ``stress-energy tensor of
Faddeev-Popov ghosts'' (or ``$b$ and $c$ fields of spins $n$ and 
$1-n$'') on the Riemann surface $X\simeq\Gamma\bk\Omega$.  

Thus we arrive at the main result of the paper.

\begin{theorem}\label{main-thm}
Let $g$ and $n$ be integers, $g>1$, $n>1$, 
and let ${\det}'\Delta_{n}$ and $S$ be the 
functions on Schottky space $\schottky_{g}$
defined above.  Let 
$p: \univschottky_{g}\rightarrow \schottky_{g}$
be the universal Schottky curve, let 
$T_{V}\univschottky_{g}$ be the vertical tangent bundle,
and let $\varphi_{1},\dots,\varphi_{d_{n}}$
be the family of global holomorphic sections of
$p_{*}(T_{V}\univschottky_{g}^{-n})$ 
(the``natural basis'' for $n$-differentials)
defined in Section \ref{zeromodes} below, forming a basis
for each fibre.  
For $t\in\schottky_{g}$ let
$[N_{n}]_{jk}(t)=\langle\varphi_{j}(t),\varphi_{k}(t)\rangle$,
where the inner product is induced from the 
hyperbolic metric on the compact Riemann
surface $X_{t}\simeq\Gamma_{t}\bk\Omega_{t}$. 
Then there exists a 
holomorphic function $F(n):\schottky_{g}\to\C$
such that
\begin{equation} \label{n-factorization}
\frac{{\det}'\Delta_{n}}{\det N_{n}}=
 c_{g,n}
 \exp\{-\frac{6n^{2}-6n+1}{12\pi}S\}
  \abs{F(n)}^{2},
\end{equation}
where $c_{g,n}$
is a constant depending only on $g$ and $n$.
The function $F(n)$ is given by the following absolutely 
convergent product:
\begin{equation} \label{F-n-formula}
F(n)=(1-q_{L_{1}})^{2}\dots(1-q_{L_{1}}^{n-1})^{2}(1-q_{L_{2}}^{n-1}) \prod_{\{\gamma\}}\prod_{m=0}^{\infty}(1-q_{\gamma}^{n+m}),
\end{equation}
where $q_{\gamma}$ is the multiplier of $\gamma\in\Gamma_{t}$,
$\{\gamma\}$ runs over all distinct primitive conjugacy classes 
in the marked normalized Schottky group $\Gamma_{t}$,
excluding the identity, and $L_{1},\dotsc, L_{g}$ are the free generators
fixing the marking of $\Gamma_{t}$.
\end{theorem}

See section \ref{kleinian-groups} for the definitions of $q_{\gamma}$ and
primitive $\gamma$, and for the normalization of the marked Schottky
group.
For $n>1$ and $g>1$, we have ${\det}'\Delta_{n} =C_{g,n}Z(n)$, 
where $Z(s)$ is the Selberg zeta function 
\eqref{selberg-zeta} and $C_{g,n}$ is 
a constant depending only on $g$ and $n$ (\cite{DP86:DetLaplacian,Sarnak87:DetLaplacian}), so that 
Theorem \ref{main-thm} gives a factorization of $Z(n)$ for 
integers $n>1$, 
considered as functions on $\schottky_{g}$. As in the case of 
Zograf's formula, the function $F(n)$ defines an isometry between 
the line bundle $\lambda_{n}$ over $\moduli_{g}$ equipped with the 
Quillen metric, and the holomorphic line bundle over
$\moduli_{g}$ determined by the
Hermitian metric $\exp\lbrace\frac{6n^{2}-6n+1}{12\pi}S\rbrace$. 
Theorem \ref{main-thm}, together with \eqref{kahler-potential},
immediately implies the local families index theorem \eqref{local-index},
of which it may be considered the ``$\del\delbar$ antiderivative''.

Heuristically, the function $F(n)$ on $\schottky_{g}$ can be interpreted
as a holomorphic determinant ${\det}'\delbar_{n}(t)$ of the family of 
$\delbar_{n}$-operators on Riemann surfaces $X_{t},\,t\in\schottky_{g}$, in accordance with arguments in \cite{Knizhnik89:Multiloop}.  We note in passing
that the functions $F(1)$ and $F(2)$ enter the ``Polyakov measure for 
the $D=26$ theory of closed bosonic strings'' \cite{BK86:HoloAn,Knizhnik89:Multiloop,DHoker:Strings}.

The content of the paper is the following.
In Section \ref{basics} we collect the facts we 
will need on Kleinian groups, Green's functions,
Teichm\"uller and Schottky spaces, and the 
classical Liouville action.  In Section \ref{series}
we express the Green's function of $\delbar_{n}$
in terms of Poincar\'e series, thus completing the 
outline given in \cite{Mar87:CFTonRS}.  Section
\ref{zeromodes} describes our choice of a natural,
holomorphically varying basis of
$H^{0}(X_{t},\canclass^{n}_{X_{t}})$.  Finally in Section
\ref{proof} we prove Theorem \ref{main-thm}.  For $n=1$,
our proof is essentially the argument of \cite{Zograf97:FPreprint},
which establishes Theorem \ref{zograf-thm} for those 
Schottky groups with exponent of convergence $\delta<1$.
(For the first part of Theorem \ref{zograf-thm} when
$\delta\geq 1$, we refer to \cite{Zograf89:FPaper}.)

The results of this paper may be extended to the case
where the $n$-differentials on $X$ are twisted by a 
character of the Schottky group, or equivalently, a 
unitary character of $\pi_{1}(X)$, generalizing Kronecker's
second limit formula.  In this case, comparison
with known bosonization results yields a product formula
for theta functions in genus $g>1$, generalizing the 
Jacobi triple product formula when $g=1$.  We intend to 
return to this in a sequel to this paper.

\smallskip\noindent
\textbf{Acknowledgments} We are grateful to Peter Zograf for 
sharing his insights with us. We are also thankful to Lee-Peng Teo 
for useful discussions at the early stages of this work. A.M.~would 
like to thank the organizers of the Simons Workshop in Mathematics and Physics in Stony Brook during Summers 2003 
and 2004, when this work was completed, for support. The work of
L.T. ~was partially supported by the NSF grant DMS-0204628.

%----------------------------------------
\section{Necessary basic facts}
\label{basics}
%----------------------------------------

Here we fix notation, and recall the basic definitions and known results we 
will need.

%----------------------------------------
\subsection{Kleinian groups}\label{kleinian-groups}
\cite{Bers75:AutoformsSchottky}
%----------------------------------------
By definition, a \define{Kleinian group} 
is a discrete subgroup $\Gamma$ of the group of M\"obius transformations 
$\mathrm{PSL}(2,\C)$
which acts properly discontinuously on some non-empty
open subset of the Riemann sphere $\widehat{\C}=\C\cup\{\infty\}$. 
The largest such subset $\Omega\subset\widehat{\C}$
is called the \define{ordinary set} of $\Gamma$ and its complement 
is called the \define{limit set} of 
$\Gamma$. 

For integers $n$ and $m$, an 
\define{automorphic form of type $(n,m)$} for $\Gamma$ 
is a function $f:\Omega\to\widehat{\C}$ such that
\[
f(z)=f(\gamma z)\,\gamma'(z)^{n}
\overline{\gamma'(z)}^{m}\quad
\text{for all}\quad z\in\Omega,\;\gamma\in\Gamma.
\]
We write the space of smooth forms of type
$(n,m)$ as $\autoform^{n,m}\(\Omega,\Gamma\)$
(abbreviating $\autoform^{n,0}=\autoform^{n}$),
and the space of holomorphic forms
of type $(n,0)$ as $\holoautoform^{n}\(\Omega,\Gamma\)$.
A \define{function group} is a Kleinian group
which leaves some connected 
component $\Omega_{0}\subseteq\Omega$
invariant, and a \define{uniformization} of 
a Riemann surface $X$ is a function group
$\Gamma$ with invariant component
$\Omega_{0}\subseteq\Omega$ such 
that $X\simeq \Gamma\backslash\Omega_{0}$.
Since $\Omega_{0}$ is invariant, we can
define the restrictions
$\autoform^{n,m}\(\Omega_{0},\Gamma\)$
and $\holoautoform^{n}\(\Omega_{0},\Gamma\)$.

The \define{exponent of convergence} of a Kleinian
group $\Gamma$ is the infimum of $\delta\in\R$ such that 
the series $\sum_{\gamma\in\Gamma}\abs{\gamma'(z)}^{\delta}$
converges for all $z\in\Omega$.  For all Kleinian groups,
$\delta<2$.

A Kleinian group $\Gamma$ 
is called a \define{Fuchsian group} if it leaves some Euclidean disc 
invariant; we will assume the disc has been
conjugated to the upper half-plane $\U=\{z=x+iy\in\C : y>0\}$,
so that $\Gamma\subset\mathrm{PSL}(2,\R)$.

A Kleinian group $\Gamma$ is called a \define{Schottky group} 
if it is generated by $L_{1},\dotsc,L_{g}$ satisfying the following
condition:  there exist $2g$ smooth Jordan curves $C_{r}$,
$r=\pm 1,\dotsc,\pm g$, which form the oriented boundary
of a domain $D\subset\widehat\C$, such that 
$L_{r}C_{r}=-C_{-r}$, $r=1,\dotsc,g$ (the negative sign
indicating opposite orientation).  The domain $D$ is a 
fundamental region for $\Gamma$.  A Schottky group is a 
function group, and a free group on generators
$L_{1},\dotsc,L_{g}$.  Each nontrivial element $\gamma$
of $\Gamma$ is \define{loxodromic}:  there exists a 
unique number
$q_{\gamma}\in\C$ (the \define{multiplier})
such that $0<\tabs{q_{\gamma}}<1$ and $\gamma$ is conjugate
in $\mathrm{PSL}(2,\C)$ to $z\mapsto q_{\gamma}z$, 
that is, 
\[
\frac{\gamma z-a_{\gamma}}{\gamma z -b_{\gamma}}
=q_{\gamma}\frac{z-a_{\gamma}}{z -b_{\gamma}}
\]
for some $a_{\gamma},\,b_{\gamma}\in\widehat\C$
(respectively, the \define{attracting} and \define{repelling
fixed points}).  A \define{marked Schottky group} is a Schottky 
group together with an ordered set of free generators 
$L_{1},\dots, L_{g}$; it is \define{normalized}
if $a_{L_{1}}=0$, $b_{L_{1}}=\infty$,
and $a_{L_{2}}=1$.

It will be convenient to define 
$L_{-r}:=L_{r}^{-1}$, so that 
$L_{r}C_{r}=-C_{-r}$ is true for all 
$r\in\{\pm 1,\dotsc,\pm g\}$.  We abbreviate
$a_{r}:=a_{L_{r}}$, $b_{r}:=b_{L_{r}}$
and $q_{r}:=q_{L_{r}}$.  Denote by 
$D_{r}$ the connected component of 
$\widehat\C-C_{r}$ 
containing $b_{r}$,
for $r=\pm1,\dotsc,\pm g$, so that $-C_{r}$
is the oriented boundary of $D_{r}$ and
$L^{s}_{r}(D)\subseteq D_{-r}$ for $s>0$.
Since $\Gamma$ is free, every nontrivial $\gamma\in\Gamma$
has a unique expression as a reduced word,
$\gamma=L_{r_{1}}^{s_{1}}\dotsb L_{r_{m}}^{s_{m}}$,
for some $r_{j}\in\{\pm1,\dotsc,\pm g\}$,
$s_{j}>0$, $j=1,\dotsc,m$, where $|r_{j}|\neq |r_{j+1}|$
for $j=1,\dotsc,m-1$.

We collect some facts we will need about the action 
of a Schottky group on $\widehat{\C}$ below.
\begin{lemma} \label{combinatorics} 
Let $\Gamma$ be a marked Schottky group. 
With notation as above, 
the following statements hold.
\begin{itemize}
\item[(i)] For all $r\neq j$ and $s>0$,
$L_{r}^{s}\(D_{j}\)\subset D_{-r}$.
\item[(ii)]  Let
$\gamma=L_{r_{1}}^{s_{1}}\dotsb L_{r_{m}}^{s_{m}}\in\Gamma$ be 
a reduced word.
Then $a_{\gamma}\in D_{-r_{1}}$ and
$b_{\gamma}\in D_{r_{m}}$.
\item[(iii)] Let
$\gamma=L_{r_{1}}^{s_{1}}\dotsb L_{r_{m}}^{s_{m}}\in\Gamma$ be a
reduced word.
Then
\[
\gamma^{-1}\(a_{r}\)\in
\begin{cases}
D_{r_{m}} & \text{if $\gamma\neq L_{r}^{s}$  for all $s>0$},\\
D_{-r}=D_{-r_{m}} & \text{if $\gamma=L_{r}^{s}$  for some $s>0$}.
\end{cases}
\]
\end{itemize}
\end{lemma} 
\begin{proof}
Part (i) trivially follows from definitions. For part (ii), we observe that
$\gamma(D)\subseteq D_{-r_{1}}$, which immediately follows from part (i) using induction on $m$,
\[ \gamma(D)= L_{r_{1}}^{s_{1}}(L_{r_{2}}^{s_{2}}\dotsb L_{r_{m}}^{s_{m}}(D))
\subseteq L_{r_{1}}^{s_{1}}(D_{-r_{2}}) 
\subseteq D_{-r_{1}}.
\]
This shows that $a_{\gamma}\in D_{-r_{1}}$.
For $b_{\gamma}$, just
note that $b_{\gamma}=a_{\gamma^{-1}}$.
Part (iii) is also proved 
by induction on $m$.  For $m=1$, if $r_{1}\neq r$,
then $\gamma^{-1}(a_{r})\in\gamma^{-1}(D_{-r})=L^{s_{1}}_{-r_{1}}(D_{-r})
\subseteq D_{r_{1}}$, while if $r_{1}=r$, then
$a_{r}$ is fixed by $\gamma$ and 
$\gamma^{-1}\(a_{r}\)=a_{r}\in D_{-r}$.  Now assume
for $m-1$ and suppose $\gamma\neq L_{r}^{s}$ 
for all $s>0$.  Then
\[
\gamma^{-1}(a_{r})
=L_{-r_{m}}^{s_{m}}(
 (L_{r_{1}}^{s_{1}}\dotsb L_{r_{m-1}}^{s_{m-1}})^{-1}
 (a_{r})) \\
\in L_{-r_{m}}^{s_{m}}(D_{\pm r_{m-1}})
 \subseteq D_{r_{m}}.
\]
\end{proof}
For future use, we mention that an element $\gamma$
of a group $\Gamma$ is called \define{primitive}
if $\gamma\neq\gamma_{0}^{s}$ for all 
$\gamma_{0}\in\Gamma$ and integers $s>1$.

%----------------------------------------
\subsection{The operators $\delbar_{n}$ and $\Delta_{n}$}
 \label{greensfns}
%----------------------------------------
We follow \cite{TakZog87:Localindexthm}.  Let $X$ be a 
compact Riemann surface of genus $g>1$.  X carries a 
unique hyperbolic metric (a Hermitian metric of constant 
curvature $-1$), written locally as $\rho(z)\tabs{\ud z}^2$.
Let $\canclass_{X}=T^{\ast}X$ be the holomorphic cotangent 
bundle of $X$, i.e., the canonical class, 
and for any integers $n$ and $m$, let 
$\diffform^{p,q}(X,\canclass_{X}^n\otimes\overline{\canclass}_{X}^m)$
be the vector space of smooth differential forms of type $(p,q)$ 
on $X$ with values in the line bundle $\canclass_{X}^n\otimes\overline{\canclass}_{X}^m$.
An \define{$(n,m)$-differential} (or \define{$n$-differential} 
when $m=0$) is an element of $\autoform^{n,m}(X)
=\diffform^{0,0}(X, \canclass_{X}^n\otimes\overline{\canclass}_{X}^m)$
(or $\autoform^{n}(X)$ when $m=0$), written locally as 
$\varphi(z)(\ud z)^n(\ud\zbar)^m$.  Note that 
we may identify
$\diffform^{p,q}(X,\canclass_{X}^{n}\otimes\overline{\canclass}_{X}^{m})
\simeq\autoform^{n+p,m+q}(X)$.
When $X\simeq\Gamma\bk\Omega_{0}$ for some function 
group $\Gamma$ and invariant component $\Omega_{0}$, 
we identify
$\autoform^{n,m}(X)\simeq\autoform^{n,m}\(\Omega_{0},\Gamma\)$. 
In what follows we will make implicit identifications of this kind 
without further comment.

The hyperbolic metric on $X$ induces a Hermitian metric 
\begin{equation} \label{metric-n-diff}
\langle\varphi,\psi\rangle
=\iint\limits_D \varphi\overline{\psi}\rho^{1-n-m}\,\ud^{2}z,
\end{equation}
on  
$\autoform^{n,m}(X)$, where $D$ is a fundamental region for 
$\Gamma$ in $\Omega_{0}$, and 
$\ud^{2}z=\frac{i}{2}\ud z\wedge\ud\zbar$ is the Euclidean
area form on $\Omega_{0}$.  The metric and complex structure
determine a connection 
\[
D=\del_n\oplus\delbar_n:
\diffform^{0,0}(X,\canclass_{X}^n)\to
\diffform^{1,0}(X,\canclass_{X}^n)\oplus\diffform^{0,1}(X,\canclass_{X}^n)
\] 
on the line bundle $\canclass_{X}^{n}$, given locally by
\[
\delbar_n=\frac{\del}{\del\zbar}\quad\text{and}\quad
\del_n=\rho^n\,\frac{\del}{\del z}\,\rho^{-n}.
\]
The metric determines $\delbar$-Laplacians
$\Delta_{n}=\Delta_n^{0,0} = \delbar_n^*\delbar_n$
and $\Delta_{n,1}=\Delta_n^{0,1} = \delbar_n\delbar_n^*$,
acting on vector spaces $\autoform^{n}(X)$ and $\autoform^{n,1}(X)$ 
respectively, where $\delbar_n^*=-\rho^{-1}\del_n$
is the adjoint of $\delbar_{n}$ with respect to
\eqref{metric-n-diff}.

Let $\Ltwoautoform^{n,m}(X)$ be the $L^{2}$-closure of 
$\autoform^{n,m}(X)$ with respect to the inner product 
\eqref{metric-n-diff}.
The operators $\Delta_{n}$ and $\Delta_{n,1}$ are self-adjoint and 
non-negative, and have pure 
discrete spectrum in the Hilbert spaces $\Ltwoautoform^{n}(X)$ 
and $\Ltwoautoform^{n,1}(X)$. The corresponding eigenvalues 
$0\leq\lambda_0\leq\lambda_1\leq\dotsb$
of $\Delta_{n}$ (the non-zero eigenvalues of $\Delta_{n}$ and $\Delta_{n,1}$ coincide) have finite multiplicity and accumulate only at infinity. 
The determinant of $\Delta_n$ is defined by zeta regularization: the 
elliptic operator zeta-function 
\[
\zeta_{n}(s)=\sum_{\lambda_k>0}\lambda_k^{-s},
\]
defined initially for $\re s>1$, has a meromorphic continuation to the 
entire $s$-plane \cite{MP49:Eigenfns}, and by definition 
\cite{RS71:Torsion,RS73:Complextorsion},
\[
\det\Delta_n=e^{-\zeta_{n}'(0)}.
\]
The non-zero spectrum of $\Delta_{1-n}$ is identical to that of  
$\Delta_{n,1}$ (see, e.g., \cite{TakZog87:Localindexthm}), so that 
$\det\Delta_{n}=\det\Delta_{1-n}$. Hence without loss of generality
we will usually assume $n\geq 1$.

Denote by $I_{n}$ and $P_{n}$, respectively, the identity operator 
in $\Ltwoautoform^{n}(X)$, and the orthogonal projection operator 
from $\Ltwoautoform^{n}(X)$ onto
$\holoautoform^{n}(X)=\ker\delbar_{n}=\ker\Delta_{n}$. 
The \define{Green's operators} for $\delbar_{n}$
and $\Delta_{n}$ for $n\geq 1$ are the unique operators
$K_{n}:\Ltwoautoform^{n,1}(X)\to\Ltwoautoform^{n}(X)$ and 
$G_{n}:\Ltwoautoform^{n}(X)\to\Ltwoautoform^{n}(X)$ respectively,
such that 
\begin{itemize}
\item[\textbf{GF1.}] $K_{n}\delbar_{n}=G_{n}\Delta_{n}=I_{n}-P_{n}$.
\item[\textbf{GF2.}]
$\left. K_{n}\right\rvert_{\ker\delbar_{n}^{*}}
=0$ and $\left. G_{n}\right\rvert_{\ker\Delta_{n}}
=0$.
\end{itemize}
They are related by $K_{n}=G_{n}\delbar_{n}^{*}$.  Now, let $X\simeq\Gamma\bk\Omega_{0}$ for some function group 
$\Gamma$ and invariant component $\Omega_{0}$.  The 
\define{Green's functions} for $\delbar_{n}$ and $\Delta_{n}$ 
are the unique automorphic forms in two variables $K_{n}(z,z')$ 
and $G_{n}(z,z')$ respectively, smooth for $z'\neq\gamma z$, 
$z,z'\in\Omega_{0}$ and $\gamma\in\Gamma$, satisfying
\begin{align*} 
(K_{n}\psi)(z)
& =\iint\limits_{D}K_{n}(z,z')\psi(z')~\ud^{2}z' \quad\text{for all}\quad \psi\in
\autoform^{n,1}\(\Omega_{0},\Gamma\)\\
\text{and}\quad
(G_{n}\psi)(z)
&=\iint\limits_{D}G_{n}(z,z')\psi(z')~\ud^{2}z'\quad\text{for all}\quad \psi\in
\autoform^{n}\(\Omega_{0},\Gamma\).
\end{align*}
The form $K_{n}(z,z')$ is of type $(n,0)$ in $z$ and type $(1-n,0)$ in 
$z'$, and the form $G_{n}(z,z')$ is of type $(n,0)$ in $z$ and type 
$(1-n,1)$ in $z'$.  Both forms are holomorphic in $z$.  The relation $K_{n}=G_{n}\delbar_{n}^{\ast}$ implies
\[
K_{n}(z,z')=-\(\delbar_{1-n}'\)^{*}G_{n}(z,z')
=\rho(z')^{-n}\frac{\del}{\del z'}\(\rho(z')^{n-1} G_{n}(z,z')\).
\]
\begin{remark}
Our convention differs 
from \cite{TakZog87:Localindexthm}, where the Green's function
$\widetilde{G}_{n}(z,z')$ is defined by 
$(G_{n}\psi)(z)=\langle\widetilde{G}_{n}(z,\cdot),\psi\rangle$.
The two are related by 
$G_{n}(z,z')=\rho(z')^{1-n}\widetilde{G}_{n}(z,z')$.
\end{remark}

The Green's function $Q_n(z,z')$ for $\Delta_n$ on the upper 
half plane $\U$ is uniquely determined by the following properties:
\begin{enumerate}
\item[\textbf{1.}] $Q_n(z,z')$ is smooth for $z\neq z'$;
\item[\textbf{2.}] $Q_n(\gamma z,\gamma z')
      \gamma'(z)^n\gamma'(z')^{1-n}\overline{\gamma'(z')}=Q_n(z,z')$ 
      for all $\gamma\in\text{PSL}(2,\R)$ and $z\neq z'$;
\item[\textbf{3.}] 
$Q_n(z,z')=-\frac{1}{\pi}(\im z')^{-2}\log\abs{z-z'}^2 + \bigO(1)$ as $z\to z'$;
\item[\textbf{4.}] $\Delta_n Q_n(z,z')=0$ for $z\neq z'$;
\end{enumerate}
and an additional growth condition as $z\to\del\U$ (see \cite{Hejhal:TF}).  
The terminology is justified since if $X\simeq\Gamma\bk\U$ for a 
Fuchsian group $\Gamma$, then
\begin{equation*}
G_n(z,z')=\sum_{\gamma\in\Gamma}Q_n(z,\gamma z')
                \gamma'(z')^{1-n}\overline{\gamma'(z')}.
\end{equation*}
Correspondingly, the Green's function $R_{n}(z,z')$ for 
$\delbar_{n}$ on $\U$ is
$R_n(z,z')=-(\delbar_{1-n}')^{*}Q_{n}(z,z')$, and 
from the defining properties of 
$Q_{n}(z,z')$ we derive
\begin{equation} \label{R-n}
R_n(z,z')
=\frac{1}{\pi}\cdot\frac{1}{z-z'}\left(\frac{\zbar-z'}{\zbar-z}\right)^{2n-1}.
\end{equation}

%----------------------------------------
\subsection{Teichm\"uller and Schottky spaces}
\label{variational}
\cite{Bers72:Survey,Bers75:AutoformsSchottky,Hejhal:Schottky}
%----------------------------------------
A \define{marked Riemann surface} is a compact Riemann surface 
$X$ of genus $g>1$, equipped with (up to an inner automorphism of 
$\pi_{1}(X,x_{0})$) a canonical system of generators $\alpha_{1},\dots,\alpha_{g},\beta_{1},\dots,\beta_{g}$ of 
$\pi_{1}(X,x_{0})$, i.e., a system with the single relation
$\alpha_{1}\beta_{1}\alpha_{1}^{-1}\beta_{1}^{-1} \cdots\alpha_{g}\beta_{g}\alpha_{g}^{-1}\beta_{g}^{-1}=1$.
Marked Riemann surfaces will be denoted by 
$[X]=(X; \alpha_{1},\dots,\alpha_{g},\beta_{1},\dots,\beta_{g})$.
Let $\teich_{g}$ be the Teichm\"{u}ller space of marked Riemann 
surfaces of genus $g>1$.

For a marked Riemann surface $[X]$, let $\mathcal{N}$ be the smallest 
normal subgroup in $\pi_{1}(X,x_{0})$ containing 
$\alpha_{1},\dots,\alpha_{g}$.
By the classical retrosection theorem, there exists a 
Schottky group $\Gamma\simeq\pi_{1}(X,x_{0})/\mathcal{N}$
with ordinary set $\Omega$ such that 
$X\simeq\Gamma\bk\Omega$.
The group $\Gamma$ is unique if we require it to be 
normalized; we will always assume that $\Gamma$ 
is normalized  and marked by generators $L_{1},\dots,L_{g}$ 
corresponding to the cosets 
$\beta_{1}\mathcal{N},\dotsc,\beta_{g}\mathcal{N}$. The correspondence
\[
[X]\mapsto(a_{3},\dots,a_{g},b_{2},\dots,b_{g},q_{1},\dots,q_{g})
\]  
defines a complex-analytic map $\Psi:\teich_{g}\to\C^{3g-3}$.
Its image $\schottky_{g}=\Psi(\teich_{g})$ is a domain in $\C^{3g-3}$, 
called the \define{Schottky space}, and $\Psi$ is a covering
map onto $\schottky_{g}$. The correspondence 
$t\mapsto\Gamma_{t}\bk\Omega_{t}$
defines a complex-analytic covering map 
$\schottky_{g}\rightarrow\moduli_{g}$.  

Equivalently, the Schottky space $\schottky_{g}$ may be defined as 
the set of marked, normalized Schottky groups of rank $g>1$, 
with a complex structure described as follows. 
For every $t\in\schottky_{g}$ let $X_{t}\simeq \Gamma_{t}\bk\Omega_{t}$
be the corresponding Riemann surface, and let $\delbar_{2}(t)$ 
and $\Delta_{-1}^{0,1}(t)$ be as defined in section
\ref{greensfns}, for the surface $X_{t}$.  Then the  holomorphic 
tangent space $T_{t}\schottky_{g}$ is naturally isomorphic to 
$\holoautoform^{-1,1}(\Omega_{t},\Gamma_{t})=\ker\Delta_{-1}^{0,1}(t)
\subset\autoform^{-1,1}(\Omega_{t},\Gamma_{t})$ ---
the space of harmonic Beltrami differentials --- 
while the holomorphic cotangent space 
$T^{\ast}_{t}\schottky_{g}$ is naturally isomorphic to
$\holoautoform^{2}(\Omega_{t},\Gamma_{t})
=\ker\delbar_{2}(t)
\subset\autoform^{2}(\Omega_{t},\Gamma_{t})$ ---
the space of holomorphic quadratic differentials. For
$\mu\in\holoautoform^{-1,1}(\Omega_{t},\Gamma_{t})$
and $q\in\holoautoform^{2}(\Omega_{t},\Gamma_{t})$,
the pairing is given by
\[
(\mu,q)=\iint\limits_{D_{t}}\mu\,q\,\ud^{2}z,
\]   
where $D_{t}$ is a fundamental region for $\Gamma_{t}$. The inner 
product \eqref{metric-n-diff} 
on harmonic $(-1,1)$-differentials
defines a Hermitian metric on the Schottky space $\schottky_{g}$. 
This metric is K\"{a}hler, and coincides with the projection onto 
$\schottky_{g}$ of the Weil-Petersson metric on $\teich_{g}$ (see \cite{Ahlfors61:RemarksonTeich}). We will 
call it the Weil-Petersson metric on $\schottky_{g}$ and will denote its 
symplectic form by $\omega_{\scriptscriptstyle{WP}}$.  

In this definition of $\schottky_{g}$, one defines complex 
coordinates for a neighbourhood of $t\in\schottky_{g}$, 
called Bers coordinates, as follows.  Given
$\mu\in\holoautoform^{-1,1}(\Omega_{t},\Gamma_{t})$
satisfying
$\norm{\mu}_\infty=\sup_{z\in\Omega_{t}}\abs{\mu(z)}<1$,
there exists a unique homeomorphism
$f^\mu:\widehat\C\to\widehat\C$ fixing
$0,1,\infty$ and satisfying the Beltrami equation 
\[
\frac{\del f^\mu}{\del \zbar}=\mu\frac{\del f^\mu}{\del z}.
\]
Set $\Gamma^\mu=f^\mu\circ\Gamma\circ (f^\mu)^{-1}$, $\Omega^{\mu}=f^{\mu}(\Omega)$, and
$X^\mu=\Gamma^\mu\bk\Omega^{\mu}$. Choosing a basis 
$\mu_1,\dotsc,\mu_{3g-3}$ for 
$\holoautoform^{-1,1}(\Omega_{t},\Gamma_{t})$ gives 
$\mu=\vep_{1}\mu_{1} +\dots +\vep_{3g-3}\mu_{3g-3}$, 
where $\vep_{i}\in\C$. The correspondence  
$\vep=(\vep_{1},\dots,\vep_{3g-3})\mapsto\Psi([X^{\mu}])$ 
introduces complex coordinates in a neighborhood of 
$t\in\schottky_{g}$; the corresponding complex structure
agrees with that given by the first definition, considering
$\schottky_{g}$ as a domain in $\C^{3g-3}$.  In terms of Bers 
coordinates,
\[
\omega_{\scriptscriptstyle{WP}}
\(\frac{\del}{\del\vep_{k}},\frac{\del}{\del\bar{\vep}_{l}}\)
=\frac{i}{2}\langle \mu_{k},\mu_{l}\rangle
\quad\text{at}~t\in\schottky_{g}.
\]

The Schottky universal curve is a fibration
$p:\univschottky_{g}\to\schottky_{g}$ with fibre 
$\pi^{-1}(t)=X_{t}\simeq\Gamma_{t}\bk\Omega_{t}$ 
for $t\in\schottky_{g}$.  Let 
$T_V\univschottky_g\to\univschottky_g$ be  the holomorphic vertical 
tangent bundle --- the  holomorphic line bundle over 
$\univschottky_g$
consisting of vectors in the holomorphic tangent space 
$T\univschottky_g$ that are tangent to the fibres $X_{t}=\pi^{-1}(t)$.  
A  family $\varphi^{\varepsilon}$ of $(n,m)$-differentials on Riemann 
surfaces $X^{\vep\mu}$ is defined as a smooth section of the line 
bundle
\[
(T_V\univschottky_g)^{-n}\otimes(\overline{T_V\univschottky_g})^{-m}
\to\univschottky_g.
\]
The hyperbolic metric $\rho$ gives rise to a family of $(1,1)$-differentials 
and defines a natural Hermitian metric on the line bundle $T_V\univschottky_g\rightarrow\univschottky_g$, whose restriction to 
each fibre coincides with the hyperbolic metric. It also defines a Hermitian 
metric in the bundle $(T_V\univschottky_g)^{-n}\to\univschottky_g$, and 
in the direct image bundle 
$\Lambda_{n}=p_{*}((T_V\univschottky_g)^{-n})\to\schottky_{g}$. 
The fibre of $\Lambda_{n}$ over $t\in\schottky_{g}$ is the vector space $\holoautoform^{n}(\Omega_{t},\Gamma_{t})$, and the corresponding 
Hermitian metric is given by \eqref{metric-n-diff}.

The pullback of an $(n,m)$-differential $\varphi^{\varepsilon}$ 
over $X^{\varepsilon\mu}$ is an $(n,m)$-differential over 
$X=X^{0}$, defined by 
\begin{equation*}
f_*^{\varepsilon\mu}(\varphi^\varepsilon)=\varphi^\varepsilon\circ f^{\varepsilon\mu}
\,(f^{\varepsilon\mu}_{z})^n(\overline{f^{\varepsilon\mu}_{z}})^m,
\end{equation*}
where $f^{\varepsilon\mu}:\widehat\C\to\widehat\C$ is the 
corresponding solution of Beltrami equation.  The Lie derivatives of the 
family $\varphi^\varepsilon$ in the directions $\mu$ and 
$\overline{\mu}$, where 
$\mu\in\holoautoform^{-1,1}(\Omega_{t},\Gamma_{t})
\simeq T_{t}\schottky_{g}$ 
and $t=\Psi([X])$, are defined by
\begin{align*}
\delta_\mu\varphi
&=\frac{\del}{\del\varepsilon}\bigg\vert_{\varepsilon=0}
 f_*^{\varepsilon\mu}(\varphi^\varepsilon) \in\autoform^{n,m}(X)
\\
\text{and}\quad
\bar{\delta}_\mu\varphi
&=\frac{\del}{\del\overline{\varepsilon}}\bigg\vert_{\varepsilon=0}
 f_*^{\varepsilon\mu}(\varphi^\varepsilon) \in\autoform^{n,m}(X).
\end{align*}
Every smooth function $\varphi$ on $\schottky_{g}$ is naturally identified 
with a family of $(0,0)$-differentials, constant along the fibres of $p$, 
which we will continue to denote by $\varphi$. In this case the Lie derivative coincides with the usual directional derivative,
\begin{align*}
\var_\mu\varphi=\del\varphi(\mu)
\quad \text{and}\quad
\varbar_\mu\varphi=\delbar\varphi(\mu),
\end{align*}
where $\del$ and $\delbar$ are the $(1,0)$ and $(0,1)$ components,
respectively, of the deRham differential $\ud$ on the complex manifold  $\schottky_{g}$.
Similarly, for a family of linear operators 
$A^\varepsilon:\autoform^{k,l}(X^{\vep\mu})\to\autoform^{m,n}(X^{\vep\mu})$ 
we define the Lie derivatives by 
\begin{align*}
\delta_\mu A
&=\frac{\del}{\del\varepsilon}\bigg\vert_{\varepsilon=0}
 \left(f_*^{\varepsilon\mu}A^\varepsilon(f_*^{\varepsilon\mu})^{-1}\right)
\\
\text{and}\quad
\bar{\delta}_\mu A
&=\frac{\del}{\del\overline{\varepsilon}}\bigg\vert_{\varepsilon=0}
 \left(f_*^{\varepsilon\mu}A^\varepsilon(f_*^{\varepsilon\mu})^{-1}\right),
\end{align*}
so that
\[
\delta_{\mu}(A(\vp)) =\delta_{\mu}A(\vp)+ A(\delta_{\mu}\vp)\quad\text{and}\quad\bar{\delta}_{\mu}(A(\vp)) =\bar{\delta}_{\mu}A(\vp)+ A(\bar{\delta}_{\mu}\vp).
\]

Now we present some variational formulas we will need. 
For $\mu\in\holoautoform^{-1,1}\(\Omega,\Gamma\)$ 
define
\begin{equation*}
F_\mu=\frac{\del}{\del\varepsilon}f^{\varepsilon\mu}\bigg\vert_{\varepsilon=0}
\quad\text{and}\quad
\Phi_\mu=\frac{\del}{\del\bar{\varepsilon}}
          f^{\varepsilon\mu}\bigg\vert_{\varepsilon=0}.
\end{equation*}
Then \cite{Ahlfors61:RemarksonTeich}
\[
\frac{\del F_\mu}{\del \bar{z}} =\mu \quad\text{and}\quad \Phi_\mu=0,
\] 
and $\displaystyle{\chi_{\mu}[\gamma]=\frac{F_{\mu}\circ\gamma}{\gamma'} -F_{\mu}}$ is a polynomial of order $\leq 2$ every $\gamma\in\Gamma$.
(For groups other than
Schottky, function $\Phi_{\mu}$ is holomorphic on $\Omega$ but not
necessarily zero.)  Note that the normalization of $f^{\vep\mu}$ implies
that $F_{\mu}(0)=F_{\mu}(1)=F_{\mu}(\infty)=0$, and hence
$\chi_{\mu}[L_{1}](0)=0$, $\chi_{\mu}[L_{1}](\infty)=0$, and
$\chi_{\mu}[L_{2}](1)=0$.  (Here $F_{\mu}(\infty)=0$ means
$F_{\mu}(z)=\littleo(\tabs{z}^{2})$ as $z\to\infty$, and similarly for 
$\chi_{\mu}[L_{1}]$.)  Another classical result of Ahlfors \cite{Ahlfors61:RemarksonTeich} is that for the family $\rho$ of $(1,1)$-differentials given by the hyperbolic metric, 
\[
\delta_\mu\rho = 0\quad\text{and}\quad \bar{\delta}_{\mu}\rho=0.
\]
>From this one finds (see, e.g., \cite{TakZog87:Localindexthm}),
\[
\delta_\mu\delbar_n=-\mu\del_n\quad\text{and}\quad \delta_\mu\del_n=0,
\]
and hence
\[
\delta_\mu\Delta_n=\rho^{-1}\mu\del_{n+1}\del_{n}.
\]
If $\varphi$ is a smooth family of holomorphic automorphic forms of
type $(n,0)$, then differentiating $\delbar_{n}\varphi=0$ one gets
\begin{equation} \label{var-phi}
\delbar_{n}(\var_\mu\varphi)=\mu\del_n\varphi\quad\text{and}\quad
\delbar_{n}(\bar{\delta}_{\mu}\vp)=0,
\end{equation} 
where the last equation follows from $\bar{\delta}_{\mu}\delbar_{n}=0$.
Finally, for $t\in\schottky_{g}$ let $\gamma_{t}\in\Gamma_{t}$ be a group 
element corresponding to a fixed element $[\gamma]$ under the isomorphism 
$\Gamma_{t}\simeq\pi_{1}(X,x_{0})/\mathcal{N}$. Then the multipliers 
$q_{\gamma_{t}}$ give rise to a holomorphic function $q_{\gamma}:\schottky_{g}\to\C$. 
Identifying 
$T_{t}^{*}\schottky_{g}\simeq\holoautoform^{2}(\Omega_{t},\Gamma_{t})$, we have
(see e.g. \cite{Zograf89:FPaper})
\begin{equation}\label{dq}
\del q_{\gamma}=-\frac{q_{\gamma}}{\pi}
\sum_{\sigma\in\langle\gamma\rangle\backslash\Gamma}
\frac{\(a_{\gamma}-b_{\gamma}\)^{2}}
     {\(\sigma z-a_{\gamma}\)^{2}\(\sigma z-b_{\gamma}\)^{2}}
   \sigma'(z)^{2},
\end{equation}
where the sum runs over the set of left cosets in $\Gamma$
of the cyclic subgroup generated by $\gamma$. 

%----------------------------------------
\subsection{Classical Liouville action}
\label{liouville-action-section}
\cite{TakZog87:Liouvilleaction,TakTeo03:Liouvilleaction}
%----------------------------------------
The Schottky space $\schottky_{g}$ is a domain of holomorphy
\cite{Hejhal:Schottky}, so that the Weil-Petersson metric on 
$\schottky_{g}$ has a globally defined K\"{a}hler potential. Here 
we present the potential for the Weil-Petersson metric constructed 
in \cite{TakZog87:Liouvilleaction}. It is given by the 
``classical Liouville action'' --- the critical value of the 
``Liouville action functional'' on the family of Riemann surfaces 
parameterized by the Schottky space $\schottky_{g}$ ---
and has the additional property of establishing a relation between 
Fuchsian and Schottky uniformizations. 
 
Namely for $t\in\schottky_{g}$ set $X=X_{t}$; for convenience, 
we omit the subscript $t$ here and write $X\simeq\Gamma\bk\Omega$, 
etc. 
Let $\rho(z)\tabs{\ud z}^{2}$ be the hyperbolic metric on $\Omega$,
pulled back from the hyperbolic metric on $X\simeq\Gamma\bk\Omega$. 
Let $D$ be a fundamental region for the marked Schottky group 
$\Gamma$ (see section 2.1).
Set
\begin{align*}
S &= \iint\limits_{D}
  \(\abs{\frac{\del\log\rho}{\del z}}^{2}+\rho\)\,\ud^{2} z \\
&\quad+ \frac{i}{2}\sum_{k=2}^{g}
 \oint_{C_{k}}\(\log\rho-\frac{1}{2}\log\abs{L_{k}'}^{2}\)
  \(\frac{L_{k}''}{L_{k}'}~\ud z -
   \frac{\overline{L_{k}''}}{\overline{L_{k}'}}~\ud\zbar\)
   \\
&\quad+4\pi\sum_{k=2}^{g}\log\abs{c\(L_{k}\)}^{2},
\end{align*}
where for $\gamma=
\bigl(\begin{smallmatrix}
a&b\\c&d\\
\end{smallmatrix}\bigr)\in\Gamma$,
we denote $c(\gamma)=c$. This definition does not depend on a 
particular choice of the fundamental region $D$.
The values $S_{t}$ for $t\in\schottky_{g}$ define a smooth 
function $S:\schottky_{g}\to\R$, called the \emph{classical Liouville action} 
(see \cite{TakZog87:Liouvilleaction} for motivation and details, and \cite{TakTeo03:Liouvilleaction} for a cohomological interpretation). 
The function $S$ is invariant with respect to transformations of 
$\schottky_{g}$ corresponding to permutations of the generators 
of the marked Schottky group \cite{Zograf89:FPaper}.
For a holomorphic function $f$ with $f'\neq 0$, the
\define{Schwarzian derivative} of $f$ is
\begin{equation}\label{schwarz}
\schwarz(f)=\(\frac{f''}{f'}\)' -\frac{1}{2}\(\frac{f''}{f'}\)^{2}.
\end{equation}
For $X\simeq \Gamma\bk\Omega$ let 
$J:\U\to\Omega$ be the universal covering of $\Omega$ and set 
\[
\vartheta =2\schwarz(J^{-1}).
\]
Though the mapping $J$ is not one-to-one, it follows from the properties of $J$ and $\schwarz$ that $\vartheta$ is a well-defined element 
of $\holoautoform^{2}(\Omega,\Gamma)$ \cite{TakZog87:Liouvilleaction}. Correspondingly, the smooth family $\vartheta_{t}$ of holomorphic quadratic differentials on $X_{t}$ gives rise to a $(1,0)$-form $\vartheta$ on 
$\schottky_{g}$. 
\begin{proposition} \label{var-liouvilleaction} 
The function $S:\schottky_{g}\to\R$ has the followng properties.
\begin{itemize}
\item[(i)] $\del S=\vartheta$;
\item[(ii)] $\del\delbar S =2i\,\omega_{\scriptscriptstyle{WP}}$.
\end{itemize}
\end{proposition}
\begin{proof}
See \cite{TakZog87:Liouvilleaction} (and \cite{TakTeo03:Liouvilleaction} for generalization to Kleinian groups of class A).
\end{proof}

%----------------------------------------
\section{Poincar\'e series and the Green's function of $\delbar_{n}$}
\label{series}
%----------------------------------------

Let $X\simeq\Gamma\bk\Omega_{0}$ for some function group 
$\Gamma$ and invariant component $\Omega_{0}$, and let 
$n$ be a positive integer.  In this 
section we define a meromorphic Poincar\'e series 
$\widehat{K}_n(z,z')$ and a smooth kernel $K_n^0(z,z')$ associated 
with the subspace 
$\holoautoform^{n}(\Omega_{0},\Gamma)=\ker\delbar_{n}$,
such that for $n>1$ the Green's function $K_n\(z,z'\)$ of $\delbar_n$ 
is given by $K_n=\widehat{K}_n + K_n^0$.  (There is a slight modification
when $n=1$.)  This completes the outline 
sketched in \cite{Mar87:CFTonRS}.

For convenience, assume that $\infty$ is in the limit set of $\Gamma$.
For $n>1$, fix points $A_1,\dotsc,A_{2n-1}$ in the limit set of 
$\Gamma$, such that 
\begin{equation}\label{Aj-condition}
\forall j~\exists~\text{at most}~ n-1~\text{distinct}
~k~\text{such that}~A_{k}=A_{j}.
\end{equation}
If $n=1$, fix a single point $A_{1}$ in the ordinary set of 
$\Gamma$.
Then for $n\geq 1$ and $z,z'\in\Omega_{0}$ with $z'\neq\gamma z$ for all 
$\gamma\in\Gamma$, define \cite{Bers67:IneqforFGKleiniangps}
\begin{equation} \label{series-bers}
\widehat{K}_n(z,z')
=\frac{1}{\pi}\sum_{\gamma\in\Gamma}\frac{1}{\gamma z-z'}
 \biggl(\prod_{j=1}^{2n-1}\frac{z'-A_j}{\gamma z-A_j}\biggr)
 \gamma'(z)^n,
\end{equation}
with the natural conventions if $A_{j}=\infty$ for one or more $j$. 
\begin{lemma} Let $\Gamma$ and $\widehat{K}_{n}$ be defined as 
above.
\begin{itemize}
\item[(i)] Suppose $n>1$. For $z,z'\in\Omega_{0}$ with $z'\neq\gamma z$ 
for all $\gamma\in\Gamma$, the series $\widehat{K}_n(z,z')$ converges 
absolutely and uniformly on compact subsets. It defines a meromorphic
function on $\Omega_{0}\times\Omega_{0}$ with only simple poles, 
at $z'=\gamma z$, $\gamma\in\Gamma$. 
\item[(ii)] Suppose that $\Gamma$ has exponent of convergence $\delta<1$.  
Then for $z,z'\in\Omega_{0}$ with 
$z'\neq\gamma z$ and $z\neq\gamma A_{1}$ for all
$\gamma\in\Gamma$, the series $\widehat{K}_1(z,z')$ converges 
absolutely and uniformly on compact subsets. It defines a meromorphic
function on $\Omega_{0}\times\Omega_{0}$ with only simple poles, 
at $z'=\gamma z$ and $z=\gamma A_{1}$, $\gamma\in\Gamma$.
\end{itemize}
\end{lemma}
\begin{proof} Since $L^{1}$-convergence of holomorphic functions 
implies uniform convergence on compact sets, for (i) it is sufficient to 
show 
\[
\iint\limits_{\Omega_{0}}\frac{1}{|z-z'|}
\prod_{j=1}^{2n-1}\frac{1}{|z-A_j|}\,
\rho(z)^{1-n/2}\ud^{2}z 
<\infty,
\]
where $\rho(z)\tabs{\ud z}^{2}$ is the hyperbolic metric on $\Omega_{0}$.
This was proved in \cite{Bers71:EichlerIntwSingularities} using Ahlfors'
estimates for $\rho$, under the assumption that 
$A_{1},\dotsc,A_{2n-1}$ are distinct points in the limit set.
Exactly the same proof works when some of the $A_{j}$ coincide,
provided they satisfy condition \eqref{Aj-condition}.  Because
$A_{1}$ is in the ordinary set for $n=1$, (ii) follows immediately 
from the definition of $\delta$.
\end{proof}

Let $\Pi_{2n-2}$ be the vector space of polynomials of degree 
$\leq 2n-2$, considered as a right $\Gamma$-module with the 
$\gamma\in\Gamma$ acting on $p\in\Pi_{2n-2}$ by
\[
\gamma_{*}p = p\circ \gamma\cdot(\gamma')^{1-n},
\]
and denote by $Z^{1}(\Gamma,\Pi_{2n-2})$ the vector space of  
$1$-cocycles for the group $\Gamma$ with coefficients in $\Pi_{2n-2}$ 
--- the \define{Eichler cocycles} \cite{Bers67:IneqforFGKleiniangps}. 
Explicitly, a cocycle is a map $\chi: \Gamma\to\Pi_{2n-2}$ satisfying
\[
\chi[\gamma_{1}\gamma_{2}] ={\gamma_{2}}_{*}\chi[\gamma_{1}] +\chi[\gamma_{2}]\;\;\,\text{for all}\;\;\,\gamma_{1},\gamma_{2}\in\Gamma.
\] 
A direct computation shows that for any $\gamma\in\Gamma$,
\begin{align*}
\widehat{K}_n(\gamma z,z')\gamma'(z)^{n}
&=\widehat{K}_n(z,z')\\
\widehat{K}_n(z,\gamma z')\gamma'(z')^{1-n}
&=\widehat{K}_n(z,z')
 + \chi_{\widehat{K}}[\gamma](z,z'),
\end{align*} 
where $\chi_{\widehat{K}}(z,\,\cdot\,)\in Z^{1}(\Gamma,\Pi_{2n-2})$ 
for every $z\in\Omega_{0}$, and 
$\chi_{\widehat{K}}[\gamma](\,\cdot\,,z')
\in\holoautoform^{n}(\Omega_{0},\Gamma)$ 
for every $\gamma\in\Gamma$ and $z'\in\C$.

Now, let $\varphi_{1},\dotsc,\varphi_{d}$ be a basis for 
$\holoautoform^{n}(\Omega_{0},\Gamma)$, where $d=(2n-1)(g-1)$ 
($n>1$), or $d=g$ ($n=1$).  Define \define{potentials} $F_{k}$ (there should be no confusion with $F_{\mu}$ defined in section 
\ref{variational})  of the 
automorphic forms $\varphi_{k}$ by 
\cite{Bers67:IneqforFGKleiniangps,Bers71:EichlerIntwSingularities}
\begin{equation}\label{potential-definition}
\begin{split}
F_k(z)
&=
-\frac{1}{\pi}\iint\limits_{\Omega_0}
 \frac{\rho(\zeta)^{1-n}\overline{\varphi}_k(\zeta)}{\zeta-z}
 \prod_{j=1}^{2n-1}\frac{z-A_j}{\zeta-A_j}\,
 \ud^2\zeta
\\
&=
-\iint\limits_{D_{0}}\rho(\zeta)^{1-n}\overline\varphi_{k}(\zeta)
\widehat{K}_{n}(\zeta,z)\,\ud^{2}z\\
&=
-\langle\widehat{K}_{n}(\,\cdot\,,z),\varphi_{k}\rangle,
\end{split}
\end{equation}
where $\rho(\zeta)$ is the hyperbolic metric on $\Omega_0$. Note that
though $\widehat{K}_{n}(\,\cdot\,,z)$ is not in $\Ltwoautoform^{n}(\Omega_{0},\Gamma)$, 
the inner product given by \eqref{metric-n-diff} is still well-defined.  
The function $F_k$ on $\Omega_{0}$ has the property 
\begin{equation} \label{potential}
\frac{\del F_k}{\del \bar{z}}=\rho^{1-n}\bar{\varphi}_k.
\end{equation}
Let $[N_n]_{jk}=\langle\varphi_j,\varphi_k\rangle$ be the Gram matrix 
of the basis $\varphi_{1},\dots,\varphi_{d}$ with respect to the inner 
product \eqref{metric-n-diff}, and let $N_{n}^{jk}=[N_{n}^{-1}]_{jk}$ be 
the inverse matrix. For $z,z'\in\Omega_{0}$ set
\begin{equation}\label{zero-mode-kernel}
K_n^0(z,z')=\sum_{j=1}^d\sum_{k=1}^d
N_n^{kj}\varphi_j(z) F_k(z').
\end{equation}
It follows from \eqref{potential} that
\begin{equation} \label{projection}
\frac{\del K_n^0}{\del \bar{z}'}(z,z') =P_n(z,z')
\end{equation}
is the integral kernel of the orthogonal projection 
$P_{n}:\Ltwoautoform^{n}(\Omega_{0},\Gamma)
\to\holoautoform^{n}(\Omega_{0},\Gamma)$.  
For any $\gamma\in\Gamma$ we have
\begin{align*}
K_n^0(\gamma z,z')\gamma'(z)^{n}
&=K_n^0(z,z')\\
K_n^0(z,\gamma z')\gamma'(z')^{1-n}
&=K_n^0(z,z')-\sum_{j=1}^{d}\sum_{k=1}^{d}N_n^{kj}\varphi_{j}(z)
 \langle\chi_{\widehat{K}}[\gamma](\,\cdot\,,z'),\varphi_{k}\rangle\\
&=K_n^0(z,z')
 - \chi_{\widehat{K}}[\gamma](z,z'),
\end{align*} 
since $\chi_{\widehat{K}}[\gamma](\,\cdot\,,z')
\in\holoautoform^{n}(\Omega_{0},\Gamma)$.
Hence $\widehat{K}_{n}+K_n^0$ is an automorphic form of type 
$(n,0)$ in $z$ and type $(1-n,0)$ in $z'$.
\begin{proposition} 
Let $\Gamma$, $\widehat{K}_{n}$ and $K_{n}^{0}$
be defined as above, and let $K_{n}$ be the Green's function for 
$\delbar_{n}$ on $\Gamma\bk\Omega_{0}$ defined in section
\ref{greensfns}.  Then:
\begin{itemize}
\item[(i)] for $n>1$ and $z,z'\in\Omega_{0}$,
\[
K_{n}(z,z')=\widehat{K}_{n}(z,z')+K_n^0(z,z');
\]
\item[(ii)] if $\delta<1$, then for $z,z'\in\Omega_{0}$,
\[
K_{1}(z,z')-K_{1}(z,A_{1})=\widehat{K}_{1}(z,z')+K_1^0(z,z').
\]
\end{itemize}
\end{proposition}
\begin{proof}
First we verify condition \textbf{GF1}, i.e., show that for any $\varphi\in\autoform^{n}\(\Omega_{0},\Gamma\)$,
\[
\iint\limits_{D_{0}}
(\widehat{K}_n+K_n^0)(z,z')(\delbar_{n}\varphi)(z')\ud^{2}z'
=\varphi(z)-(P_n\varphi)(z),
\]
where $D_{0}$ is a fundamental region for $\Gamma$ in $\Omega_{0}$.  
We have
\[
\iint\limits_{D_{0}}
 (\widehat{K}_n+K_n^0)(z,z')(\delbar_{n}\varphi)(z')\ud^{2}z'
= I_{1} - I_{2},
\]
where
\begin{align*}
I_{1} &= \lim_{\varepsilon\to 0}\iint
 \limits_{D_{0}\setminus\{|z'-z|\leq\vep\}}
 \delbarprime_{n}\((\widehat{K}_n+K_n^0)(z,z')\varphi(z')\)\ud^{2}z', \\
I_{2} & = \lim_{\varepsilon\to 0}\iint
 \limits_{D_{0}\setminus\{|z'-z|\leq\vep\}}
 \delbarprime_{n}\((\widehat{K}_n+K_n^0)(z,z')\)\varphi(z')\ud^{2}z'.
\end{align*}
By Stokes' theorem, $I_{1}$ is a sum of an integral over the boundary of 
$D_{0}$, which vanishes since $(\widehat{K}_n+K_n^0)(z,z')\varphi(z')$
is a $(1,0)$-differential in $z'$, and a boundary term around the singularity 
$z'=z$, so that
\[
I_{1} = \lim_{\varepsilon\to 0}\,
 \frac{1}{2\pi i}\oint\limits_{\abs{z'-z}=\varepsilon}
 \(\frac{\varphi(z')}{z'-z} + \bigO(1)\)\ud z'
=\varphi(z).
\] 
Since $\widehat{K}_n(z,z')$ is holomorphic in $z'$ for $z'\neq z$, using \eqref{projection} we get 
$I_{2}=(P_{n}\varphi)(z)$.

Since condition \textbf{GF2} is vacuous for $n>1$, the above establishes
(i) in that case.  When $n=1$, the above argument shows that the operators $K_{1}$
and $\widehat{K}_{1}+K_{1}^{0}$ agree on $\im\delbar_{1}$, that is,
\[
K_{1}(z,z')=\widehat{K}_{1}(z,z')+K_{1}^{0}(z,z')+\psi(z)
\]
for some $\psi\in\holoautoform^{1}(\Omega_{0},\Gamma)$.  Setting
$z'=A_{1}$ evaluates $\psi$ and yields (ii).
\end{proof}

\begin{remark} It follows that
\[
\frac{\del K_1}{\del z'}(z,z')
= \frac{\del \widehat{K}_1}{\del z'}(z,z')+ \frac{\del K^0_1}{\del z'}(z,z'),
\]
which is Fay's formula relating Bergmann and Schiffer kernels on a 
compact Riemann surface \cite{Fay77:Fouriercoeff}.  This was used in 
the proof of the local families index theorem \eqref{local-index} in the 
case $n=1$ given in \cite{TakZog87:Localindexthm}, and was the starting 
point for the proof of Zograf's factorization
formula \eqref{F-0} in \cite{Zograf97:FPreprint}.
\end{remark}

%----------------------------------------
\section{Natural basis for $H^{0}(\schottky_g,\Lambda_n)$}
\label{zeromodes}
%----------------------------------------

It was proved by Kra \cite{Kra84:Poincare} that the direct 
image vector bundle
\[
\Lambda_n=p_{*}((T_V\univschottky_g)^{-n})\to\schottky_g
\] 
is holomorphically trivial, i.e., there exist 
$\varphi_1,\dots,\varphi_d\in H^{0}(\schottky_g,\Lambda_n)$ 
such that for each $t\in\schottky_g$, the holomorphic $n$-differentials 
$\varphi_1(t),\dots,\varphi_d(t)$ on $X_t$ form a basis of the fibre
$\holoautoform^n(X_t)$. 
For $n=1$, the abelian differentials $\varphi_1(t),\dots,\varphi_g(t)$ 
on the Riemann surface $X_{t}$ with the classical normalization
\[
\oint_{\alpha_{k}}\varphi_{j}=\delta_{jk}
\]
form such a basis, since every $t\in\schottky_g$ uniquely determines 
the $\alpha$-cycles on the Riemann surface $X_t=\Gamma_t\bk\Omega_t$ 
(see \cite{Zograf89:FPaper}). Here we construct a natural basis of 
the global sections of $\Lambda_{n}$ for $n>1$, which reduces to 
the former when $n=1$.    

Let $\Gamma$ be normalized, marked Schottky group with 
distinguished system of generators $L_{1},\dotsc,L_{g}$.  For $n>1$, 
a cocycle $\chi\in Z^{1}(\Gamma,\Pi_{2n-2})$ is called \define{normalized} 
if 
\[
\frac{\del^{r}\chi[L_{1}]}{\del z^{r}}(z)=0,\;\; 0\leq r\leq n-2,
\quad 
\chi\left[L_{1}\right](z)=\littleo(\abs{z}^{n})\;\;
\text{as}\;\;z\to\infty,
\]
and $\chi[L_{2}](1)=0$. Every cocycle $\chi\in Z^{1}(\Gamma,\Pi_{0})=Z^{1}(\Gamma,\C)$
is called normalized by definition.  Let $\widetilde{Z}^{1}(\Gamma,\Pi_{2n-2})$ 
be the vector space of normalized Eichler cocycles.  Since any cocycle 
may be normalized by adding a coboundary  
$b\in B^{1}(\Gamma,\Pi_{2n-2})$
--- a cocycle $b[\gamma]=\gamma_{*}p -p$ for some $p\in\Pi_{2n-2}$ ---   
and every normalized $b\in B^{1}(\Gamma,\Pi_{2n-2})$ is identically 
zero, we have an isomorphism
\[
H^{1}(\Gamma,\Pi_{2n-2})
:=Z^{1}(\Gamma,\Pi_{2n-2})/B^{1}(\Gamma,\Pi_{2n-2})
\simeq\widetilde{Z}^{1}(\Gamma,\Pi_{2n-2}).
\]
Let 
$\Pi_{2n-2}^{g}=\underbrace{\Pi_{2n-2}\times\cdots\times\Pi_{2n-2}}_{g}$, 
and define
\[
\widetilde{\Pi}_{2n-2}^{g}
=\{\(p_{1},\dotsc,p_{g}\)\in\Pi_{2n-2}^{g}
:p_{1}(z)=cz^{n-1},~p_{2}(1)=0\}.
\]

Since the group $\Gamma$ is free,  the mapping from 
$\widetilde{Z}^{1}(\Gamma,\Pi_{2n-2})$ to $\widetilde{\Pi}_{2n-2}^{g}$
given by
\[
\chi\mapsto (\chi[L_{1}],\dotsc,\chi[L_{g}])
\]
is an isomorphism.  Fix a basis of $\widetilde{\Pi}_{2n-2}^{g}$; this fixes 
a basis 
\[
\xi_{1},\dotsc,\xi_{d}
\in\widetilde{Z}^{1}(\Gamma,\Pi_{2n-2})\simeq H^{1}(\Gamma,\Pi_{2n-2}).
\]  
This basis depends only on $\Gamma$ as an abstract group --- that is, $\xi_{k}[\gamma]$ depends only on the reduced word 
$L_{r_{1}}^{s_{1}}\dotsb L_{r_{m}}^{s_{m}}$
representing $\gamma$. Thus we have defined a basis of 
$H^{1}(\Gamma,\Pi_{2n-2})$ 
simultaneously for all normalized marked Schottky groups 
$\Gamma_{t},\, t\in\schottky_{g}$.

Now we define a basis for $\mathcal{H}^{n}(\Omega,\Gamma)$ 
corresponding to the basis $\xi_{1},\dots,\xi_{d}$ of
$\widetilde{Z}^{1}(\Gamma,\Pi_{2n-2})$ associated with a fixed basis of
$\widetilde{\Pi}_{2n-2}^{g}$. For this purpose we use the Bers map 
$\beta^{*}:\holoautoform^{n}(\Omega,\Gamma)
\to H^{1}(\Gamma,\Pi_{2n-2})$, where $\chi=\beta^{*}(\varphi)$ is defined by 
\[
\chi[\gamma] =F\circ\gamma\cdot(\gamma')^{1-n} -F,
\]
with $F$ a potential of the holomorphic $n$-differential $\varphi$ 
given by \eqref{potential-definition}. The potential $F$ depends on the 
points $A_{1},\dots,A_{2n-1}$ in the limit set of $\Gamma$; a 
different choice of normalization points adds a coboundary to 
$\chi$. We will always choose the normalization points to be $\underbrace{0,\dots,0}_{n-1},1,\underbrace{\infty,\dots,\infty}_{n-1}$. 
With this normalization, we get a mapping 
\[
\widetilde{\beta}^{*}:
\holoautoform^{n}(\Omega,\Gamma)
\to\widetilde{Z}^{1}(\Gamma,\Pi_{2n-2}).
\] 
Since the Bers mapping $\beta^{*}$ is injective, $\widetilde{\beta}^{*}$ 
is also; and the vector spaces  $\holoautoform^{n}(\Omega,\Gamma)$ 
and $\widetilde{Z}^{1}(\Gamma,\Pi_{2n-2})$ have the same dimension 
$d$, so $\widetilde{\beta}^{*}$ is a complex anti-linear isomorphism. 
Define a basis $\psi_{1},\dots,\psi_{d}$ of 
$\holoautoform^{n}(\Omega,\Gamma)$ by
\[ 
\widetilde{\beta}^{*}(\psi_{k})=\xi_{k},
\] 
and let $\varphi_{1},\dots,\varphi_{d}$ be the dual basis of $\holoautoform^{n}(\Omega,\Gamma)$ with respect to the inner product
\eqref{metric-n-diff}:
\[
\langle\varphi_{j},\psi_{k}\rangle=\delta_{jk}.
\]
\begin{lemma} The holomorphic $n$-differentials 
$\varphi_{1}(t),\dotsc,\varphi_{d}(t)\in\mathcal{H}^{n}(X_{t})$, 
constructed above for every point $t\in\schottky_{g}$, define global holomorphic sections $\varphi_{1},\dots,\varphi_{d}$ of the bundle $\Lambda_{n}$ over $\schottky_{g}$.
\end{lemma}  
\begin{proof}
It follows from the construction that the $\varphi_{j}$ are smooth global
sections of $\Lambda_{n}$; we must show they are holomorphic.  
Fix $t\in\schottky_{g}$ and abbreviate $\varphi_{j}(t)=\varphi_{j}$, 
$\Gamma_{t}=\Gamma$, etc.  Let 
$\mu\in\holoautoform^{-1,1}(\Omega,\Gamma)$ represent a tangent
vector at $t$. It follows from \eqref{var-phi}   
that $\delbar_{n}(\varbar_{\mu}\varphi_{j})=0$, i.e., $\varbar_{\mu}\varphi_{j}\in\mathcal{H}^{n}(\Omega,\Gamma)$.  But
by the definition of $\xi_{k}$ and Stokes' theorem,
\begin{align}\label{dual}
\delta_{jk}=\langle\varphi_{j},\psi_{k}\rangle
=\iint\limits_{D}\varphi_{i}\frac{\del F_{k}}{\del \zbar}\,\ud^{2}z
=-\frac{1}{2i}\sum_{r=1}^{g}\oint_{C_{r}}\varphi_{j}\,\xi_{k}[L_{r}]~\ud z.
\end{align}
Since $\xi^{\vep}_{k}[L_{r}^{\vep}]$ do not depend explicitly on $\vep$
and $\Phi_{\mu}=0$, we have $\bar{\delta}_{\mu}\xi_{k}[L]=0$, so
\[
0=-\frac{1}{2i}\sum_{r=1}^{g}\oint_{C_{r}}
 (\varbar_{\mu}\varphi_{j})\xi_{k}[L_{r}]\,\ud z 
 =\langle\varbar_{\mu}\varphi_{j},\psi_{k}\rangle
\]
for each $k$, and we conclude $\bar{\delta}_{\mu}\varphi_{j}=0$.
\end{proof}
\begin{remark}
It is necessary to take the dual basis $\varphi_{j}$
because the $\psi_{k}$ are \emph{not} holomorphic sections of the bundle $\Lambda_{n}\rightarrow\schottky_{g}$. This is related to the fact that the Bers mapping $\beta^{*}$ is complex anti-linear.
\end{remark}
We say that the sections $\varphi_{1},\dots,\varphi_{d}$ form a 
\define{natural basis of $H^{0}(\schottky_{g},\Lambda_{n})$ 
corresponding to the basis $\xi_{1},\dots,\xi_{d}$ of $\widetilde{Z}^1(\Gamma,\Pi_{2n-2})$ associated with a fixed basis of $\widetilde{\Pi}^{g}_{2n-2}$} (for brevity, a \define{natural basis}).
Note that for $n=1$, if we make the choice 
\[
\xi_{k}[L_{r}]=-2i\delta_{kr},
\]
we recover the classical normalized basis of abelian differentials;
we add this condition to the definition of natural basis when $n=1$.

The vector bundle $\Lambda_{n}\to\schottky_{g}$ has a Hermitian 
metric defined by the inner product \eqref{metric-n-diff} on the fibres $\holoautoform^{n}(\Omega_{t},\Gamma_{t}),\,t\in\schottky_{g}$, 
which induces a Hermitian metric $\norm{\,\cdot\,}_{n}^{2}$
on its determinant line bundle $\lambda_{n}=\wedge^{d}\Lambda_{n}$. 
The natural basis gives a global holomorphic section
$\varphi=\varphi_{1}\wedge\cdots\wedge\varphi_{d}$ of $\lambda_{n}$,
with 
\[
\norm{\varphi}_{n}^{2}=\det N_{n},
\]
where $[N_{n}]_{jk}=\langle\varphi_{j},\varphi_{k}\rangle$.  The metric
and complex structure define a connection on $\lambda_{n}$, which in the holomorphic frame given by $\varphi$ is $\ud+\del\log\det N_{n}$,
where $\ud=\del+\delbar$ is the deRham operator on $\schottky_{g}$.

When $n=1$, the connection $(1,0)$ form on $\schottky_{g}$ can be found explicitly.  By the 
Riemann bilinear relations, $N_{1}=\im\boldsymbol\tau$, and we
have Rauch's formula \cite{Rauch65:Periods} 
\[
\del\boldsymbol\tau_{jk}(\mu)
=-2i\iint\limits_{D}\varphi_{j}\varphi_{k}\mu\,\ud^{2}z
\]
for $\mu\in\holoautoform^{-1,1}(\Omega,\Gamma)$, from which we 
obtain
\begin{equation} \label{var-zeromode-1}
\del\log\det N_{1}(\mu)
=-\iint\limits_{D}\sum_{j=1}^{g}\sum_{k=1}^{g}
  N_{1}^{kj}\varphi_{j}\varphi_{k}\,\mu\,\ud^{2}z,
\end{equation}
where $N_{1}^{jk}=[N^{-1}]_{jk}$.

There is an analog of \eqref{var-zeromode-1} for the natural basis when $n>1$. 
Namely, let
\begin{equation}\label{Tn-zeromode}
T_{n}^{0}(z) 
=\left.\(n\frac{\del}{\del z'} - (1-n)\frac{\del}{\del z}\)
  K_n^0(z,z')\right|_{z'=z}, 
\end{equation}
where $K_{n}^{0}$ is given by \eqref{zero-mode-kernel}, and define
\begin{equation}\label{period-T}
\varpi_{n}[\gamma]=T_{n}^{0}\circ\gamma\cdot(\gamma')^{2} - T_{n}^{0}
\end{equation}
for each $\gamma\in\Gamma$.  Then we have the following.

\begin{proposition}\label{var-zeromode}
Let $\varphi_{1},\dotsc,\varphi_{d}$ be a natural basis of
$H^{0}(\schottky_{g},\Lambda_{n})$ as constructed above.  Fix 
$t\in\schottky_{g}$ and abbreviate $\varphi_{j}(t)=\varphi$,
$\Gamma_{t}=\Gamma$, etc.  Let $N_{n}$, $T_{n}$, $\varpi_{n}$
be defined as above, and recall the notation for the marked
normalized Schottky group $\Gamma$ fixed in section
\ref{kleinian-groups}.  Then for 
$\mu\in\holoautoform^{-1,1}(\Omega,\Gamma)\simeq T_{t}\schottky_{g}$
with potential $F_{\mu}$, 
we have
\begin{equation}\label{var-zeromode-n}
 \del\log\det N_{n}(\mu)
 =   \iint\limits_{D}T^{0}_{n}\mu\,\ud^{2}z
  + \frac{1}{2i}\sum_{r=1}^{g}
         \oint_{C_{r}}\varpi_{n}[L_{r}]F_{\mu}\,\ud z.
\end{equation}
\end{proposition}
\begin{proof} 
Using holomorphy of the family $\varphi_j$, Stokes' theorem,
$\psi_{j}=\sum_{k=1}^{d}N_{n}^{kj}\varphi_{k}$
and \eqref{var-phi}, 
we have
\begin{align*}
\del\log\det N_{n}(\mu)
&= \sum_{j=1}^d\sum_{k=1}^d N_{n}^{kj}
  \langle\var_{\mu}\varphi_j,\varphi_k\rangle
=\sum_{j=1}^d\sum_{k=1}^d N_{n}^{kj}
  \iint\limits_{D}(\var_{\mu}\varphi_j)\frac{\del F_k}{\del\zbar}\,\ud^2 z \\
&=
-\iint\limits_{D}\left.(\del_{n}K_n^0)\right\vert_{\Delta}\mu\,\ud^2 z
 - \frac{1}{2i}\sum_{r=1}^{g}\oint_{C_{r}}
       \sum_{j=1}^d (\var_{\mu}\varphi_{j})\xi_{j}[L_{r}]\,\ud z,
\end{align*}
where $\Delta$ stands for the restriction on the diagonal $z'=z$.
This implies
\[
\begin{split}
\del\log\det N_{n}(\mu)
&= \iint\limits_{D}T_{n}^{0}\mu\,\ud^2 z
 - n\iint\limits_{D}\del_{1}(\left.K_n^0\right\rvert_{\Delta})\mu\,\ud^2 z \\
&\quad - \frac{1}{2i}\sum_{r=1}^{g}\oint_{C_{r}}
       \sum_{j=1}^{d}(\var_{\mu}\varphi_{j})\xi_{j}[L_{r}]\,\ud z,
\end{split}
\]
since
$T_{n}^0=-\left.(\del_{n}K_n^0)\right\vert_{\Delta}
 + n\del_{1}(K_n^0\vert_{\Delta})$.
Using Stokes' theorem again and $\del_{-1}\mu=0$, we obtain
\[
\iint\limits_{D}\del_1(\left.K_{n}^{0}\right\rvert_{\Delta})\mu\,\ud^2 z
=\frac{1}{2i}\sum_{r=1}^g\oint_{C_r}\sum_{j=1}^d
\varphi_{j}\,\xi_{j}[L_{r}]\mu\,\ud\zbar.
\]
Hence we must show that 
\begin{equation}\label{norm-relation}
\sum_{r=1}^g\oint_{C_r}\varpi[L_{r}]F_{\mu}\,\ud z\\
=-\sum_{r=1}^g\oint_{C_r}\sum_{j=1}^d
 (\var_{\mu}\varphi_{j})\xi_{j}[L_{r}]\,\ud z
   +n\varphi_{j}\,\xi_{j}[L_{r}]\mu\,\ud\zbar.
\end{equation}
It follows from \eqref{dual} that
\[
\sum_{r=1}^{g}\oint_{C_{r}}
(\var_{\mu}\varphi_{j})\xi_{k}[L_{r}]\,\ud z
+ \varphi_{j}(\var_{\mu}\xi_{k}[L_{r}])\,\ud z 
+ \varphi_{j}\,\xi_{k}[L_{r}]\mu\,\ud\zbar
=0,
\]
and we have
\[
\begin{split}
\var_{\mu}\xi_{k}[L_{r}] 
&=\left.\frac{\del}{\del \vep}\right|_{\vep=0}\xi_{k}[L_{r}]\circ f^{\vep\mu}
\cdot(f^{\vep\mu}_{z})^{1-n}\\
&=\frac{\del\xi_{k}[L_{r}]}{\del z}F_{\mu}
 +(1-n)\xi_{k}[L_{r}]\frac{\del F_{\mu}}{\del z},
\end{split}
\]
since, by construction, $\xi_{k}^{\vep}[L_{r}^{\vep}]$ does not depend 
explicitly on $\vep$.
Using the identity
\[
\begin{split}
0
&=\oint_{C_{r}}\ud
   (\varphi_{j}\xi_{k}[L_{r}]F_{\mu})\\
&=\oint_{C_{r}}\frac{\del}{\del z}
   (\varphi_{j}\xi_{k}[L_{r}]F_{\mu})\,\ud z
+ \varphi_{j}\xi_{k}[L_{r}]\,\mu\,\ud \zbar,
\end{split}
\]
we obtain
\begin{multline*}
-\sum_{r=1}^g\oint_{C_r}
  (\var_{\mu}\varphi_{j})\xi_{k}[L_{r}]\,\ud z
  + n\varphi_{j}\,\xi_{k}[L_{r}]\mu\,\ud\zbar\\
=\sum_{r=1}^g\oint_{C_r}
\(n\varphi_{j}\frac{\del \xi_{k}[L_{r}]}{\del z}
  -(1-n)\frac{\del\varphi_j}{\del z}\xi_{k}[L_{r}]\)F_{\mu}\,\ud z.
\end{multline*}
Now, a straightforward computation shows that
\begin{equation*}
\varpi_{n}[\gamma]
=\sum_{j=1}^{d}\(n\varphi_{j}\frac{\del \xi_{j}[\gamma]}{\del z}
-(1-n)\frac{\del \varphi_{j}}{\del z}\,\xi_{j}[\gamma]\),
\end{equation*}
which establishes \eqref{norm-relation} and completes the proof.
\end{proof}
To show the agreement of \eqref{var-zeromode-n} with 
\eqref{var-zeromode-1} when $n=1$, it suffices to observe that 
for this case, the properties of the potential $F_{k}$ of the basis 
element $\varphi_{k}$ imply that
\[
F_{k}(z)
=\overline{\int_{A_{1}}^{z}\varphi_{k}(\zeta)\,\ud\zeta}
-\int_{A_{1}}^{z}\varphi_{k}(\zeta)\,\ud\zeta.
\]

%----------------------------------------
\section{Proof of Theorems 1 and 2}
\label{proof}
%----------------------------------------

Since the functions $\det\Delta_{n}$, $\det N_{n}$ and $S$ 
on the Schottky space $\schottky_{g}$ are real-valued and the
function $F(n)$ on $\schottky_{g}$ is holomorphic, to prove 
Theorems \ref{zograf-thm} and \ref{main-thm} it sufficient to show 
that
\begin{equation} \label{main}
\del\log\det\Delta_{n} -\del\log F(n)
=\del\log\det N_{n} -\frac{6n^{2}-6n+1}{12\pi }\del S
\end{equation}
at all points in $\schottky_{g}$.
The $(1,0)$ forms on $\schottky_{g}$ appearing on the right hand side 
of \eqref{main} are given by Propositions \ref{var-liouvilleaction} and 
\ref{var-zeromode}. Here we complete the proof by computing the 
$(1,0)$ forms on the left hand side.

%----------------------------------------
\subsection{Computation of $\del\log\det\Delta_{n}$}
\label{vardetdeltasection}
%----------------------------------------
Let $X$ be a compact Riemann surface, with 
$X\simeq\Gamma\bk\Omega_{0}$ for some function group 
$\Gamma$ with invariant component $\Omega_{0}$, and let 
$\rho(z)\tabs{\ud z}^{2}$ be the hyperbolic metric on $\Omega_{0}$.
Define 
\begin{equation}\label{Tn}
T_{n}(z)
=\lim_{z'\to z}\(n\frac{\del}{\del z'}-(1-n)\frac{\del}{\del z}\)
 \(K_{n}(z,z')-\frac{1}{\pi}\frac{1}{z-z'}\),
\end{equation}
where $K_{n}$ is the Green's function for $\delbar_{n}$ on 
$\Gamma\bk\Omega_{0}$ defined in section \ref{greensfns}. 
When $\Omega_{0}=\U$, we will denote 
$T_{n}=T_{n}^{\mathrm{Fuchs}}$.  It easy to see that 
$T_{n}^{\mathrm{Fuchs}}\in\autoform^{2}(\U,\Gamma)$. 
Indeed, it follows from \eqref{R-n} that 
\[
\(n\frac{\del}{\del z'}-(1-n)\frac{\del}{\del z}\)R_{n}(z,z')
=\frac{1}{\pi}\frac{1}{(z-z')^{2}}  + \bigO(z-z')
\]
as $z'\to z$, so that 
\[
T_{n}^{\mathrm{Fuchs}}(z)
=\lim_{z'\to z}\(n\frac{\del}{\del z'}-(1-n)\frac{\del}{\del z}\)
 (K_{n}(z,z')- R_{n}(z,z')).
\]
It follows from property \textbf{2} in section \ref{greensfns} that 
$(K_{n}-R_{n})\rvert_{\Delta}$ is a $(1,0)$ form, and the identity
\begin{equation}\label{kernel-identity} 
T_{n}^{\mathrm{Fuchs}}
=-\left.\bigl(\del_{n}(K_{n}-R_{n})\bigr)\right\rvert_{\Delta}
 +n\del_{1}\bigl(\left.(K_{n}-R_{n})\right\rvert_{\Delta}\bigr)
\end{equation} 
proves the claim.  Here $\Delta$ stands for the restriction
on the diagonal $z' = z$.  

\begin{lemma} \label{projconn} Let $X\simeq\Gamma\bk\Omega_{0}$ for a function group $\Gamma$ with invariant component $\Omega_{0}$, let $J:\U\to\Omega_{0}$ be the holomorphic covering map of $\Omega_{0}$ 
by $\U$, and let $T_{n}$ and $T_{n}^{\mathrm{Fuchs}}$ be defined
as above. Then on $\Omega_{0}$,
\[
T_{n} = T^{\mathrm{Fuchs}}_{n}\circ J^{-1}\cdot((J^{-1})')^{2} 
+ \frac{6n^{2}-6n+1}{6\pi}\schwarz(J^{-1}),
\]
where $\schwarz$ denotes the Schwarzian derivative \eqref{schwarz}.
In particular, $T_{n}\in\autoform^{2}(\Omega_{0},\Gamma)$.
\end{lemma}
\begin{proof}
Note that while $J^{-1}$ is multiple-valued, the right side is a 
well-defined element of $\autoform^{2}(\Omega_{0},\Gamma)$.
The equality follows from the identity
\begin{equation*}
\lim_{z'\to z}\(n\frac{\del}{\del z'}-(1-n)\frac{\del}{\del z}\)\!
 \(\frac{J'(z)^{n}J'(z')^{1-n}}{J(z)-J(z')}-\frac{1}{z-z'}\)\\
=\frac{6n^{2}-6n+1}{6}\schwarz(J),
 \end{equation*}
which is verified by direct computation. This is the classical result when $n=1$. 
\end{proof}
\begin{remark} In conformal field theory, this result is known as the
statement that ``$b$-$c$ system with spins $n$ and $1-n$ has central
charge $6n^{2}-6n +1$'' (see, e.g., \cite{DHoker:Strings} and references therein).
\end{remark}
\begin{proposition}\label{var-detdelta} 
Let $\det\Delta_{n}$ be the
function on the Schottky space $\schottky_{g}$ defined in 
section \ref{greensfns}, and let $\vartheta$ be the $(1,0)$ form
on $\schottky_{g}$ defined in section \ref{liouville-action-section}.
For each $t\in\schottky_{g}$, abbreviate 
$T_{n}=T_{n}(t)$, $\Omega=\Omega_{t}$, $\Gamma=\Gamma_{t}$, etc.
Then for $\mu\in\holoautoform^{-1,1}(\Omega,\Gamma)\simeq T_{t}\schottky_{g}$, 
\[
\del\log\det\Delta_{n}(\mu)
=\iint\limits_{D}T_{n}\mu\,\ud ^{2}z 
-\frac{6n^{2}-6n+1}{12\pi}\vartheta(\mu).
\]
\end{proposition}
\begin{proof}
Set $\displaystyle{\mu^{\mathrm{Fuchs}}=\mu\circ J\,\frac{\overline{J'}}{J'}}$. It follows from Lemma \ref{projconn} that it is sufficient to prove

\[
\del\log\det\Delta_{n}(\mu)
=\iint\limits_{D} T^{\mathrm{Fuchs}}_{n}\,\mu^{\mathrm{Fuchs}}\,\ud^{2}z,
\]
where $D\subset\U$ is a fundamental region for a Fuchsian group uniformizing the Riemann surface $X\simeq\Gamma\backslash\Omega$. 
Using the identity \eqref{kernel-identity} and 
$\del_{-1}\mu=0$, this reduces to the statement
\[
\del\log\det\Delta_{n}(\mu)
= -\iint\limits_{D}
  \left.\bigl(\del_{n}(K_{n}-R_{n})\bigr)\right\rvert_{\Delta}\mu\,\ud^{2}z,
\]
which is Theorem 1 in \cite{TakZog87:Localindexthm}.  
\end{proof}

%----------------------------------------
\subsection{Computation of $\del\log F(n)$}
\label{varFsection}
%----------------------------------------
Let $\Gamma$ be a marked, normalized Schottky group.  For 
positive integer $n$ define
\begin{equation}\label{def-F0}
F_{0}(n)=\prod_{\{\gamma\}}
 \prod_{m=0}^{\infty}\(1-q_{\gamma}^{n+m}\),
\end{equation}
where $\{\gamma\}$ runs over all distinct primitive conjugacy classes 
in $\Gamma$, omitting the identity, and $q_{\gamma}$ is the multiplier 
of $\gamma$ --- see section \ref{kleinian-groups}.  The product 
converges absolutely if and only if the series 
$\sum_{\{\gamma\}}\sum_{m=0}^{\infty}\tabs{q_{\gamma}}^{m+n}$ 
converges.  One shows that this series converges provided that the 
multiplier series $\sum_{[\gamma]}\tabs{q_{\gamma}}^{n}$ converges, 
where $[\gamma]$ runs over all distinct conjugacy classes (not necessarily primitive) in $\Gamma$.  By a theorem of B\"user 
\cite{BuserJ96:MultseriesSchottky}, for a Schottky group $\Gamma$ with 
exponent of convergence $\delta$, the latter series converges provided 
$n>\delta$. It is known that $\delta<2$, hence for $n>1$ the 
product $F_{0}(n)$ converges absolutely for all Schottky groups $\Gamma$, 
and the product $F_{0}(1)$ converges absolutely provided that 
$\delta<1$.  Now we define
\begin{equation}\label{def-Fn}
F(n)=
\begin{cases}
F_{0}(1) & \text{if $n=1$},\\
\displaystyle{
 (1-q_{1})^{2}\cdots(1-q_{1}^{n-1})^{2}(1-q_{2}^{n-1})F_{0}(n)} & \text{if $n>1$}.
\end{cases}
\end{equation}

For $n\geq2$ the expression $F(n)$ defines a holomorphic function on
$\schottky_{g}$. For $n=1$ the function $F=F(1)$ is defined on the open 
subset of $\schottky_{g}$ characterized by
$\delta<1$. 

\begin{remark}
The product $\prod_{\{\gamma\}}(1-q_{\gamma}^{s})$ was briefly 
described in \cite{Bowen79:Hausdorffdimqcircles}, where it was asserted 
that with the values of $q_{\gamma}^{s}$ chosen appropriately,
the product is defined for all $\re s>\delta$ and has an analytic continuation
to the entire $s$-plane.  To our knowledge these results have not yet been 
proved.  The function $\tabs{F_{0}(n)}^{2}$ coincides with a product of 
Ruelle type zeta functions $R_{\rho}(s)$ associated to the hyperbolic
$3$-manifold $X^{3}$ defined by $\Gamma$, considered in 
\cite{Fried86:TorsiononHyperbolic}: $\tabs{F_{0}(n)}^{2}=Z_{n}(n)$, where
\[
Z_{n}(s)=\prod_{m=0}^{\infty}R_{\rho_{n+m}}(s+m),
\]
and $\rho_{n+m}$ is the representation of $\pi_{1}(X^{3})$ on 
$\mathrm{O}(2)$ taking a closed geodesic with twist parameter $\theta$ 
to a rotation of angle $(n+m)\theta$.
\end{remark}

Set
\begin{equation}\label{Tn-hol}
\widehat{T}_{n}(z)=
\lim_{z'\to z}\(n\frac{\del}{\del z'}-(1-n)\frac{\del}{\del z}\)
\(\widehat{K}_{n}(z,z')-\frac{1}{\pi}\frac{1}{z-z'}\),
\end{equation}
where $\widehat{K}_{n}$ is the Poincar\'{e} series \eqref{series-bers}. 
We have 
\begin{equation} \label{ste-decomposition}
T_{n}=\widehat{T}_{n} + T^{0}_{n},
\end{equation}
where $T^{0}_{n}$ and $T_{n}$ are defined in \eqref{Tn-zeromode} and 
\eqref{Tn} respectively.  Since 
$T_{n}\in\autoform^{2}(\Omega,\Gamma)$, we have for 
$\gamma\in\Gamma$, 
\[
\widehat{T}_{n}\circ\gamma\cdot(\gamma')^{2} - \widehat{T}_{n} 
= -\varpi_{n}[\gamma],
\]
where $\varpi_{n}[\gamma]$ is given by \eqref{period-T}.
\begin{proposition} \label{var-Fn}
Let $F(n):\schottky_{g}\to\C$ be defined by \eqref{def-F0} and
\eqref{def-Fn}.  Fix $t\in\schottky_{g}$ and abbreviate
$\Gamma_{t}=\Gamma$, etc. Let $\widehat{T}_{n}$ 
and $\varpi_{n}$ be defined by \eqref{Tn-hol} and \eqref{period-T}
respectively, corresponding to $X_{t}=X=\Gamma\bk\Omega$,
and recall the notation for the marked
normalized Schottky group $\Gamma$ fixed in section
\ref{kleinian-groups}.  For 
$\mu\in\holoautoform^{-1,1}(\Omega,\Gamma)\simeq T_{t}\schottky_{g}$
with potential $F_{\mu}$, the $(1,0)$ form $\del\log F(n)$
satisfies
\[
\del\log F(n)(\mu)
= \iint\limits_{D}\widehat{T}_{n}\mu\,\ud^{2}z
-\frac{1}{2i}\sum_{r=1}^{g}\oint_{C_{r}}\varpi[L_{r}]F_{\mu}\,\ud z.
\]
\end{proposition}
\begin{proof}
For $\gamma\in\Gamma$,
$\gamma\neq \mathrm{id}$, 
and $z\in\Omega$, we introduce
the abbreviations
\begin{align*}
A_{\gamma}(z)
&=\frac{1}{\pi}
 (nq_{\gamma}^{n-1}+(1-n)q_{\gamma}^{n})
 \frac{\gamma'(z)}{(\gamma z-z)^{2}},\\
B_{\gamma}(z)
&=\lim_{z'\to z}\frac{1}{\pi}
 \(n\frac{\del}{\del z'}  - (1-n)\frac{\del}{\del z}\)
 \frac{1}{\gamma z-z'}
 \(\prod_{j=1}^{2n-1}
  \frac{z'-A_{j}}{\gamma z-A_{j}}\)
 \gamma'(z)^{n},
\end{align*}
and split the computation into three steps.

\noindent
\emph{Step 1.} Claim that the right hand side can be written as
\begin{equation}\label{T-hol}
\iint\limits_{D}\widehat{T}_{n}\mu\,\ud^{2}z
-\frac{1}{2i}\sum_{r=1}^{g}\oint_{C_{r}}\!\!\!\varpi[L_{r}]F_{\mu}\,\ud z
= -\frac{1}{2i}
     \sum_{\substack{\gamma\in\Gamma\\\gamma\neq\mathrm{id}}}
      \,\sum_{r=1}^{g}\oint_{C_{-r}}\!\!\!\!B_{\gamma}
         \chi_{\mu}[L_{-r}]\,\ud z .
\end{equation}
We have
\[
\begin{split}
\iint\limits_{D}\widehat{T}_{n}\mu\,\ud^{2}z
&=\iint\limits_{D}\delbar(\widehat{T}_{n}F_{\mu})\,\ud^{2}z
=\frac{1}{2i}\sum_{r=1}^{g}\left(\oint_{C_{-r}}\!\!\!\!+\oint_{C_{r}}\right)
 \widehat{T}_{n}F_{\mu}\,\ud z  \\
&=-\frac{1}{2i}\sum_{r=1}^{g}\oint_{C_{r}}
\bigl((\widehat{T}_{n}-\varpi_{n}[L_{r}])(F_{\mu}+\chi_{\mu}[L_{r}])
 -\widehat{T}_{n}F_{\mu}\bigr)\,\ud z  \\
&=-\frac{1}{2i}\sum_{r=1}^{g}\oint_{C_{r}}
  \widehat{T}_{n}\circ L_{r}(L_{r}')^{2}\chi_{\mu}[L_{r}]\,\ud z
+\frac{1}{2i}\sum_{r=1}^{g}\oint_{C_{r}}
  \varpi_{n}[L_{r}]F_{\mu}\,\ud z.
\end{split}
\]
But for any Eichler cocycle, 
$\displaystyle{\chi[\gamma^{-1}]
=-\frac{\chi[\gamma]\circ\gamma^{-1}}{(\gamma^{-1})'}}$,
so we have 
\[
\oint_{C_{r}}\widehat{T}_{n}\circ L_{r}(L_{r}')^{2}\chi_{\mu}[L_{r}]\,\ud z
=\oint_{C_{-r}}\widehat{T}_{n}\chi_{\mu}[L_{-r}]\,\ud z.
\]
This, together with 
$\displaystyle{
\widehat{T}_{n}(z)
=\sum_{\gamma\in\Gamma\setminus\{\mathrm{id}\}}
B_{\gamma}(z)
}$,
converging absolutely and uniformly on compact subsets of 
$\Omega$, establishes
\eqref{T-hol}.  Note that the non-automorphy of $\widehat{T}_{n}$
necessitates the use of the integral over $C_{-r}$ rather than $C_{r}$. 

\noindent
\emph{Step 2.} Computation of $\del\log F_{0}(n)$.
Claim that
\begin{equation} \label{F-0-n}
\del\log F_{0}(n)(\mu)
=-\frac{1}{2i}
 \sum_{\substack{\gamma\in\Gamma\\\gamma\neq\mathrm{id}}}
 \,\sum_{r=1}^{g}\oint_{C_{-r}}
    A_{\gamma}\chi_{\mu}[L_{-r}]\,\ud z.
\end{equation}

Indeed,
using the expression $\displaystyle{
\log F_{0}(n)
=-\sum_{\{\gamma\}}\sum_{m=1}^{\infty}
\frac{1}{m}\frac{q_{\gamma}^{mn}}{1-q_{\gamma}^{m}}
}$ and the series \eqref{dq},
we get
\[
\begin{split}
\del\log F_{0}(n)&=\frac{1}{\pi}\sum_{\{\gamma\}}
  \sum_{\sigma\in\langle\gamma\rangle\backslash\Gamma}
  \sum_{m=1}^{\infty}
   \left[nq_{\gamma}^{m(n-1)}+(1-n)q_{\gamma}^{mn}\right]
\\
& \quad\quad\quad\cdot
   \frac{q_{\gamma}^{m}}{\(1-q_{\gamma}^{m}\)^{2}}
   \frac{\(a_{\gamma}-b_{\gamma}\)^{2}}
     {\(\sigma z-a_{\gamma}\)^{2}\(\sigma z-b_{\gamma}\)^{2}}
   \sigma'(z)^{2}
\\
&=  \frac{1}{\pi}\sum_{\{\gamma\}}
  \sum_{\sigma\in\langle\gamma\rangle\backslash\Gamma}
  \sum_{m=1}^{\infty}
   \left[nq_{\sigma^{-1}\gamma^{m}\sigma}^{n-1}
            +(1-n)q_{\sigma^{-1}\gamma^{m}\sigma}^{n}\right]
\\
&\quad\quad\quad\cdot
    \frac{1}{\(\sigma^{-1}\gamma^{m}\sigma z-z\)^{2}}
    \(\sigma^{-1}\gamma^{m}\sigma\)'(z)
\\
&= \sum_{\substack{\gamma\in\Gamma\\\gamma\neq\mathrm{id}}}
      A_{\gamma}(z),
\end{split}
\]
where we have identified 
$T_{t}^{*}\schottky_{g}\simeq\holoautoform^{2}(\Omega,\Gamma)$.
The convergence is absolute and uniform on compact subsets of 
$\Omega$.  
Since $\del\log F_{0}(n)$, unlike $\widehat{T}_{n}$, is automorphic,
applying Stokes' theorem as in step 1 gives \eqref{F-0-n}.

\noindent
\emph{Step 3.} When $n=1$, we have $\varpi[\gamma]=0$ and $A_{\gamma}(z)=B_{\gamma}(z)$, so the proposition is proved. 
For the case $n>1$ we use the assumption that the normalization points 
$A_{1},\dotsc,A_{2n-1}$ are  
$\underbrace{0,\dots,0}_{n-1},1,\underbrace{\infty,\dots,\infty}_{n-1}$ 
(see Section 4), and show that
\begin{multline}\label{normfactors}
\del\Bigl(\log\prod_{j=1}^{n-1}
 (1-q_{1}^{j})^{2}(1-q_{2}^{n-1})\Bigr)(\mu) \\
= \frac{1}{2i}
 \sum_{\substack{\gamma\in\Gamma\\
    \gamma\neq\mathrm{id}}}
 \sum_{r=1}^{g}
  \int_{C_{-r}}(A_{\gamma}-B_{\gamma})
   \chi_{\mu}\left[L_{-r}\right]\,\ud z.
\end{multline}

We first compute the right hand side of \eqref{normfactors}.  Suppose 
$\gamma\neq L_{1}^{m}$, $L_{-1}^{m}$ or $L_{2}^{m}$ for any 
$m>0$.  Direct computation verifies that
$(A_{\gamma}-B_{\gamma})(z)\chi_{\mu}[L_{-r}](z)$
is regular at $\infty$, with poles at $b_{\gamma}$, $\gamma^{-1}(0)$,
$\gamma^{-1}(1)$ and $\gamma^{-1}(\infty)$.  By part (iii) of Lemma \ref{combinatorics}, all these poles are in a single domain $D_{r_{m}}$ 
bounded by $C_{r_{m}}$ for 
$\gamma=L_{r_{1}}^{s_{1}}\dotsb_{r_{m}}^{s_{m}}$,
so that every integral in \eqref{normfactors} is zero. Thus the 
computation reduces to the cases when $\gamma = L_{1}^{m}$, 
$L_{-1}^{m}$ or $L_{2}^{m}$ for $m>0$. For $\gamma=L_{1}^{m}, m>0$,
using Lemma \ref{combinatorics} again we see that $0\in D_{-1}$ and 
$\gamma^{-1}(1), \infty\in D_{1}$. By an elementary computation, 
using the identity 
\[
\sum_{m=1}^{\infty}\frac{nq^{mn}+(1-n)q^{(n+1)m}}{(1-q^{m})^{2}}=\sum_{m=n}^{\infty}\frac{mq^{m}}{1-q^{m}},\;\;|q|<1,   
\]
and the normalization $\chi_{\mu}[L_{-1}](z)=az$, we get
\[
\frac{1}{2i}\sum_{m=1}^{\infty}
\oint_{C_{-1}}(A_{L_{1}^{m}}-B_{L_{1}^{m}})(z)
   \chi_{\mu}[L_{-1}](z)\,\ud z 
=a\sum_{j=1}^{n-1}\frac{jq_{1}^{j}}{1-q_{1}^{j}}.
\]
When $\gamma=L_{-1}^{m}, m>0$, we have $\gamma^{-1}(1), 0 \in D_{-1}$
and $\infty \in D_{1}$.  Changing $z\mapsto 1/z$ we get as before,
\[
\frac{1}{2i}\sum_{m=1}^{\infty}
\oint_{C_{-1}}(A_{L_{-1}^{m}}-B_{L_{-1}^{m}})(z)
   \chi_{\mu}[L_{-1}](z)\,\ud z
=a\sum_{j=1}^{n-1}\frac{jq_{1}^{j}}{1-q_{1}^{j}}.
\]
For $\gamma=L_{2}^{m}$ we have by Lemma \ref{combinatorics} 
that $1\in D_{-2}$ and 
$b_{2},\gamma^{-1}(0),\gamma^{-1}(\infty)\in D_{2}$. By an elementary computation, using the normalization 
$\chi_{\mu}[L_{-2}](z)=b(z-1) + c(z-1)^{2}$,
we get
\[
\frac{1}{2i}\sum_{m=1}^{\infty}
\oint_{C_{-2}}(A_{L_{2}^{m}}-B_{L_{2}^{m}})(z)
   \chi_{\mu}[L_{-2}](z)\,\ud z 
=b(n-1)\frac{q_{2}^{n-1}}{1-q_{2}^{n-1}}.
\]

To compute the left hand side of \eqref{normfactors}, we use \eqref{dq} and the identity
\[
\sum_{r=1}^{g}
 \oint_{C_{-r}}\sum_{\gamma\in\langle L\rangle\bk\Gamma}\frac{\gamma'(z)^{2}}{(\gamma z-a)^{2}(\gamma z-b)^{2}}\chi_{\mu}[L_{-r}](z)\,\ud z =\oint_{C}\frac{\chi_{\mu}[L](z)}{(z-a)^{2}(z-b)^{2}}\,\ud z,
\]
where $a=a_{L}, b=b_{L}$ and circles $C$ and $C'=-L(C)$ form the boundary for a fundamental domain of $\langle L \rangle$ in $\C\setminus\{a,b\}$. (It readily follows from Stokes' theorem and automorphy properties of the sum $\sum_{\gamma\in\langle L\rangle\bk\Gamma}$, see \cite{Kra85:Reciprocity}).
This computation establishes \eqref{normfactors} and completes 
the proof of the proposition.
\end{proof}

Theorem \ref{main-thm} now follows from \eqref{ste-decomposition} and Propositions \ref{var-liouvilleaction}, \ref{var-zeromode}, \ref{var-detdelta} 
and \ref{var-Fn} in the case $n>1$.  For $n=1$, this also gives a proof of 
Zograf's formula --- Theorem \ref{zograf-thm} --- for Schottky groups with
$\delta<1$.  For the remainder of Theorem \ref{zograf-thm} we refer 
to \cite{Zograf89:FPaper}.  
\begin{remark} Note that the functions ${\det}^{\prime}\Delta_{n}$, 
$F_{0}(n)$ and $S$ on $\schottky_{g}$ are invariant with respect to the transformations of $\schottky_{g}$ which correspond to permutations of the generators $L_{1},\dots, L_{g}$, whereas the function $\det N_{n}$ is not. 
Consequently Theorem \ref{main-thm} implies that the extra factors in the 
definition of $F(n)$ guarantee that the product $\det N_{n}\abs{F(n)}^{2}$ is invariant 
with respect to these transformations. 
This can be also verified by a direct computation.
\end{remark}

%----------------------------------------
%   Bibliography
%   Produced by BibTeX
%----------------------------------------

\providecommand{\bysame}{\leavevmode\hbox to3em{\hrulefill}\thinspace}
\providecommand{\MR}{\relax\ifhmode\unskip\space\fi MR }
\providecommand{\MRhref}[2]{%
  \href{http://www.ams.org/mathscinet-getitem?mr=#1}{#2}
}
\providecommand{\href}[2]{#2}

%----------------------------------------
%   End
%----------------------------------------

\end{document}